\documentclass[12pt]{amsart}
\usepackage{amssymb}

\newtheorem{theorem}{Theorem}[section]
\newtheorem{lemma}[theorem]{Lemma}

\newtheorem{definition}[theorem]{Definition}

\newtheorem{corollary}[theorem]{Corollary}
\newtheorem{proposition}[theorem]{Proposition}

\newtheorem{lem-def}[theorem]{Lemma-Definition}

\newcommand{\hooklongrightarrow}{\lhook\joinrel\longrightarrow}
\renewenvironment{proof}{{\bfseries Proof.}}{\qed}
\newcommand{\N}{\mathbb N}
\newcommand{\Z}{\mathbb Z}
\newcommand{\Q}{\mathbb Q}

\newcommand{\F}{\mathbb F}

\def\op{\operatorname}

\def\ab{\overline{a}}

\def\al{\alpha}

\def\ars#1{\renewcommand\arraystretch{#1}}

\def\bbr{\overline{b}}

\def\bs{\vskip.5cm}

\def\cb{\overline{c}}
\def\cc{\mathcal{C}}

\def\cfa{\left(\ri\right)_{i\in A}}

\def\chr{\op{char}}

\def\cl#1{\left[\,#1\,\right]_\mu}

\def\defn{\nn{\bf Definition. }}

\def\dg{\op{deg}}

\def\dgm{\op{deg}_\mu}

\def\dm{\Delta_\mu}

\def\diso{\lower.4ex\hbox{$\downarrow$}\raise.4ex\hbox{\mbox{\scriptsize
$\wr$}}}

\def\dta{\delta}

\def\e{\medskip}

\def\ep#1{\exp(\Pi i#1)}
\def\ep{\epsilon}

\def\erel{e_{\op{rel}}}

\def\ff{\mathcal{F}}
\def\fff{\ff(g)}
\def\ffm{\ff_\mu(g)}
\def\ffmph{\ff_{\mu,\phi}(g)}

\def\ffnq{\mathcal{F}_{\nu,Q}(g)}
\def\fin{{\op{fin}}}

\def\g{\Gamma}
\def\ga{\gamma}

\def\gal{\op{Gal}}

\def\gen#1{\big\langle\, {#1} \,\big\rangle}

\def\gg{\mathcal{G}}

\def\ggm{\mathcal{G}_\mu}
\def\ggmh{\mathcal{G}_{\muh}}

\def\gm{\g_\mu}

\def\gq{\g_\Q}
\def\gr{\operatorname{gr}}

\def\hh{{\mathcal H}}
\def\hk{\hookrightarrow}
\def\hmg{\hh(\ggm)}
\def\hmgu{\hh(\ggm^*)}
\def\hpg{\mathcal{HP}(\ggm)}

\def\hra{\hooklongrightarrow}

\def\imp{\ \Longrightarrow\ }

\def\inm{\op{in}_\mu}
\def\inmh{\op{in}_{\muh}}
\def\inn{\op{in}}

\def\irr{\op{Irr}}
\def\ism{\lower.3ex\hbox{\ars{.08}$\begin{array}{c}\,\to\\\mbox{\tiny $\sim\,$}\end{array}$}}
\def\iso{\ \lower.3ex\hbox{\ars{.08}$\begin{array}{c}\lra\\\mbox{\tiny $\sim\,$}\end{array}$}\ }

\def\ka{\kappa}
\def\kam{\kappa(\mu)}
\def\kb{\overline{K}}
\def\kh{K^h}
\def\khx{K^h[x]}
\def\km{k_\mu}

\def\kp{\op{KP}}

\def\kpi{\op{KP}_\infty}
\def\kpm{\op{KP}(\mu)}

\def\kpn{\op{KP}(\nu)}

\def\ks{K^{\op{sep}}}

\def\kx{K[x]}

\def\la{\lambda}

\def\lc{\op{lc}}

\def\lcm{\op{lc}_\mu}

\def\lfin{\ll_\fin}
\def\lg{l\raise.6ex\hbox to.2em{\hss.\hss}l}
\def\li{\ll_\infty}
\def\ll{\mathcal{L}}
\def\lra{\,\longrightarrow\,}

\def\lx{\operatorname{lex}}

\def\md#1{\; \mbox{\rm(mod }{#1})}

\def\mlv{Maclane--Vaqui\'e\ }

\def\mmu{\mid_\mu}
\def\mnu{\mid_\nu}

\def\muh{\mu^h}
\def\mul{\mu_\la}

\def\nep{\nu_\epsilon}

\def\nn{\noindent}
\def\nnn{\mathcal{N}}
\def\np{N_{\mu,\phi}}
\def\npp{N^+_{\mu,\phi}}

\def\oo{\mathcal{O}}
\def\orb{\hbox to  .3em{$\backslash$}\backslash}
\def\ord{\op{ord}}
\def\p{\mathfrak{p}}

\def\ppa{\mathcal{P}_{\al}}

\def\pset{\mathcal{P}}

\def\gqq{\mathbb{Q}\times\gq}

\def\rhc{\rho_\cc}
\def\ri{\rho_i}

\def\rj{\rho_j}

\def\sg{\sigma}

\def\sii{\ \Longleftrightarrow\ }
\def\slp{\op{sl}}

\def\smu{\sim_\mu}
\def\snu{\sim_\nu}

\def\sp{\op{Spec}}

\def\supp{\op{supp}}

\def\t{\theta}

\def\td{\,\mbox{\bf td}}
\def\tdm{\td(\mu)}
\def\th{\ttt^h}

\def\tinn{\ttt^{\op{inn}}}

\def\tmn{\ty(\mu,\nu)}

\def\ttt{\mathcal{T}}
\def\ty{\mathbf{t}}

\def\unit{\op{unit}}
\def\vb{\bar{v}}
\def\vh{v^h}

\def\wb{\overline{w}}

\newcounter{cs}
\stepcounter{cs}
\newcommand{\casos}{\begin{itemize}}
\newcommand{\fcasos}{\end{itemize}\setcounter{cs}{1}}

\newfont{\tit}{cmr12 scaled \magstep3}

\title[Polynomial factorization]{Polynomial factorization over henselian fields}

\makeatletter
\@namedef{subjclassname@2010}{%
  \textup{2010} Mathematics Subject Classification}

\author[Alberich]{Maria Alberich-Carrami$\tilde{\mbox{n}}$ana}
\address{Institut de Rob\`otica i Inform\`atica Industrial (IRI, CSIC-UPC), Institut de Mate\-mà\-tiques de la UPC-BarcelonaTech (IMTech) and Departament de Mate\-m\`a\-tiques, Universitat Polit\`ecnica de Cata\-lunya $\cdot$ BarcelonaTech, Av. Diagonal, 647, E-08028 Barcelona, Catalonia}
\email{Maria.Alberich@upc.edu}
\author[Gu\`ardia]{Jordi Gu\`ardia}
\address{Departament de Matem\`atiques, Escola Polit\`ecnica Superior d'Enginye\-ria de Vilanova i la Geltr\'u, Av. V\'\i ctor Balaguer s/n. E-08800 Vilanova i la Geltr\'u, Catalonia}
\email{jordi.guardia-rubies@upc.edu}
\author[Nart]{Enric Nart}
\address{Departament de Matem\`{a}tiques,         Universitat Aut\`{o}noma de Barcelona,         Edifici C, E-08193 Bellaterra, Barcelona, Catalonia}
\email{nart@mat.uab.cat,\ jroe@mat.uab.cat}

\author[Poteaux]{Adrien Poteaux}
\address{Univ. Lille, CNRS, Centrale Lille, UMR 9189 CRIStAL. F-59000 Lille, France}
\email{adrien.poteaux@univ-lille.fr}

\author[Ro\'e]{Joaquim Ro\'e}

\author[Weimann]{Martin Weimann}
\address{LMNO, UMR 6139, Universit\'e de Caen-Normandie, F-14032 Caen, France}
\email{martin.weimann@unicaen.fr}

\thanks{Partially supported by grants  PID2019-103849GB-I00 and PID2020-116542GB-I00  funded by MCIN/AEI/10.13039/501100011033}
\date{\today}

\begin{document}

\begin{abstract}
Given a valued field $(K,v)$ and an irreducible polynomial $g\in \kx$, we survey the ideas of Ore, Maclane, Okutsu, Montes, Vaqui\'e and Herrera-Olalla-Mahboub-Spivakovsky, leading (under certain conditions) to an algorithm to find the factorization of $g$ over a henselization of $(K,v)$.
\end{abstract}

\maketitle



\section*{Introduction}
In a pioneering work along the 1920s, \O.\!\! Ore conjectured the existence of an algorithm to compute the prime ideal decomposition of a prime number $p$ in the number field $\Q[x]/(g)$ defined by an irreducible polynomial $g\in\Q[x]$ \cite{ore1,ore2}. Ore's proposal was based in the iteration of two ``dissections":  

$\bullet$ \ Computation of Newton polygons of $g$ with respect to some valuations on $\Q[x]$.

$\bullet$ \ Factorization in certain residue fields, of residual polynomials of $g$ associated to the sides of the Newton polygons.

In the 1930s, S. Maclane solved this problem in a more general context. For a given discrete rank-one valued field $(K,v)$, he found an algorithm to compute all extensions of $v$ to the field $\kx/(g)$ defined by an irreducible polynomial $g\in\kx$ \cite{mcla,mclb}. These extensions can be identified with  certain valuations $\mu$ on $\kx$ with support $g\kx$, determined by the different irreducible factors of $g$  in $K_v[x]$, where $K_v$ is the completion  of $K$ at $v$. For each such $\mu$, Maclane constructed a chain of augmentations of valuations on $\kx$ getting arbitrarily close to it:
$$
\mu_0\;<\;\mu_1\;<\;\cdots\;<\;\mu_n\;<\;\cdots\;<\;\mu
$$
In these augmentations, some \emph{key polynomials} for the valuations $\mu_n$  are involved. This procedure can be reinterpreted as a polynomial factorization algorithm in $K_v[x]$. If a valuation $\mu_n$ is sufficiently close to $\mu$, then its key polynomial is an approximation to the irreducible factor of $g$ in $K_v[x]$ intrinsically associated to $\mu$.

Motivated by the computation of integral bases in finite extensions of local fields, K. Okutsu constructed similar approximations without using valuations on $\kx$, nor key polynomials \cite{Ok,okutsu}. 

Still in the discrete rank-one case, J. Montes developed in 1999  
certain residual polynomial operators which led to the design of a practical algorithm following the exact pattern that Ore had foreseen \cite{montes,bordeaux,HN,PW2}.
This algorithm is known as the \emph{OM-algorithm}, named after Ore, Maclane, Okutsu and Montes.

Montes' ideas led to the computation of integral bases too \cite{bases,Bauch,hayden}.
More generally, the OM-algorithm is very efficient in the resolution of many arithmetic-geometric tasks in number fields and function fields of algebraic curves \cite{gen,newapp,PW1}.

Maclane's theory was generalized to arbitrary valued fields, independently by M. Vaqui\'e \cite{Vaq0,Vaq,Vaq2,Vaq3} and F.-J. Herrera, M.-A. Olalla, W. Mahboub and M. Spivakovsky \cite{hos,hmos}.
In this general frame, \emph{limit augmentations} and the corres\-ponding \emph{limit key polynomials} appear as a new feature.

In analogy with the case of dimensions 0,1, this general theory should lead to the  development of efficient algorithms for the resolution of arithmetic-geometric tasks involving valuations of function fields of algebraic varieties of higher dimension.  Nevertheless, the extension of \cite{mclb} to a polynomial factorization algorithm over arbitrary henselian fields is still an open problem.

A prototype of general OM-algorithm was described in \cite{hmos}.
The main obstacle for this procedure to become a real algorithm is the existence of limit augmentations.

In this paper, we present an OM-algorithm for arbitrary valued fields. Let $\kh$ be a henselization of $(K,v)$. At the input of an irreducible $g\in\kx$, if the algorithm terminates, then it outputs

(a) \ Approximations to each of the irreducible factors of $g$ over $\khx$.

(b) \ All extensions of $v$ to the field $\kx/(g)$, together with a computation of their ramification indices and residual degrees.

In Sections 1--3, we review the necessary background on valuations on $\kx$, their graded algebras and Maclane--Vaqui\'e chains.

In Section 4, we discuss Newton polygons and extend Ore's dissections to this completely general setting. If $v$ has rank one, then $K$ is dense in $\kh$ and the content of this section can be  easily  deduced from Montes' original arguments in the discrete case. However, for $v$ of arbitrary rank, the key polynomials for a valuation $\mu$ on $\kx$ extending $v$ need not be irreducible over $\khx$.   
Thus, the description of the unique extension of $\mu$ to $\khx$  is more subtle and the proof of the main result (Theorem \ref{conjNewton}) is more involved.

In Section 5, we present our general OM-algorithm and we prove that, if it terminates, it solves the above mentioned problems (a) and (b).
The only obstacle for this algorithm to terminate is the existence of infinite sequences of \emph{refinement steps}. We show that there are exactly three different situations where these infinite refinements occur and we exhibit concrete examples of each one.

As an application, in Section 6 we present two OM-based algorithms where termination is guaranteed. They generalize similar constructions by Poteaux-Weimann in the discrete rank-one case \cite{PW2}.

Let $p\ge0$ be the  residual characterictic of $v$.
In Section \ref{subsecTest}, we give a general irreducibility test over $\khx$,  for all square-free polynomials $g\in\kx$ such that $p\nmid\deg(g)$. In Section \ref{subsecFactor}, restricted to the case where $v$ is (non-necessarily discrete) of rank one, we give a polynomial factorization algorithm over $\khx$, working for square-free polynomials $g\in\kx$, assuming that either $p=0$ or $\deg(g)<p$. \e

\nn{{\bf Acknowledgement. }We warmly thank Josnei Novacoski for sharing his insights about several problems in Section 4.\e

\nn{{\bf Notation. }For any field $\mathbb{K}$, we shall denote by $\irr(\mathbb{K})$ the set of all monic, irreducible polynomials in $\mathbb{K}[x]$.

\section{Commensurable extensions of a valuation to the polynomial ring}\label{secComm}


Let $(K,v)$ be a valued field, with valuation ring $\oo$ and residue class field $k$. Let $\g=v(K^*)$ be the value group and denote by $\gq=\g\otimes\Q$ the divisible hull of $\g$. 
In the sequel, we write $\gq\infty$ instead of $\gq\cup\{\infty\}$.

The equivalence classes of \emph{commensurable} extensions of $v$ to the polynomial ring $\kx$ are parametrized by the set $\ttt=\ttt(K,\gq)$  of all $\gq$-valued valuations on $\kx$, 
$$
\mu\colon \kx\lra \gq\infty,
$$
whose restriction to $K$ is $v$. 
The \emph{support} of $\mu$ is the prime ideal
$$\p=\supp(\mu):=\mu^{-1}(\infty)\in\sp(\kx).$$

The valuation $\mu$ induces a valuation $\overline{\mu}$ on the field $L$ of fractions of $\kx/\p$. That is, $L=K(x)$ if $\p=0$, or $L=\kx/\p$ if $\p=g\kx$ for some $g\in\irr(K)$. 

The residue field $\km$ of $\mu$ is, by definition, the residue field of $\overline{\mu}$.

The \emph{value group} of $\mu$ is the subgroup $\gm\subset\gq$ generated by $\mu\left(\kx\setminus\p\right)$. By definition,  
$\mu/v$ commensurable means that the quotient $\gm/\g$ is a torsion group.

We say that $\mu$ is \emph{residually transcendental} if the extension $\km/k$ is transcendental. In this case, its transcendence degree is necessarily equal to one \cite{Kuhl}.


The set $\ttt$ admits a partial ordering. 
For $\mu,\nu\in \ttt$ we say that $\mu\le\nu$ if $$\mu(f)\le \nu(f), \quad\forall\,f\in\kx.$$
As usual, we write $\mu<\nu$ to indicate that $\mu\le\nu$ and $\mu\ne\nu$.

This poset $\ttt$ has the structure of a tree. By this, we simply mean that all intervals $(-\infty,\mu\,]:=\left\{\rho\in\ttt\mid \rho\le\mu\right\}$
are totally ordered \cite[Thm. 3.9]{defless}. 

A node $\mu\in\ttt$ is a \emph{leaf} if it  is a maximal element with respect to the ordering $\le$. Otherwise, we say that $\mu$ is an \emph{inner node}.

We distinguish two kinds of leaves:  \emph{finite} and  \emph{infinite}. We denote 
$$
\ttt=\tinn\,\sqcup\,\lfin\,\sqcup\,\li
$$
the subsets of inner nodes, finite leaves, and infinite leaves, respectively. 

For all $\mu\in\ttt$, the subset to which  it belongs can be characterized as follows:

\begin{itemize}
\item $\mu\in\tinn$ if and only if $\mu$ is residually transcendental.	
\item $\mu\in\lfin$ if and only if $\supp(\mu)\ne0$.	
\item $\mu\in\li$ if and only if $\supp(\mu)=0$ and $\km/k$ is algebraic.	
\end{itemize}

The infinite leaves of $\ttt$ are \emph{valuation-algebraic}  in the terminology of Kuhlmann \cite{Kuhl}. They play no role in the polynomial factorization problem.

Let us fix an algebraic closure $\kb$ of $K$, and an extension $\vb$ of our base valuation $v$ to $\kb$. This determines a henselization $(\kh,\vh)$ of $(K,v)$. If $\ks$ is the separable closure of $K$ in $\kb$, the field $\kh\subset\ks$ is the fixed field of the decomposition group 
$$
D_{\vb}=\left\{\sg\in \gal(\ks/K)\mid \vb\circ\sg=\vb\right\}.
$$
The valuation $\vh$ is the restriction of $\vb$ to $\kh$ and it has a unique extension to $\kb$. 

\begin{theorem} \cite[Thm. A]{Rig}\label{Rig}
Let $\th=\ttt(\kh,\gq)$ be the tree of commensurable extensions of $\vh$ to $\khx$. Restriction of valuations  induces an isomorphism of posets:
$$
\th\lra\ttt,\qquad \nu\ \longmapsto \ \nu_{\mid \kx},
$$	
preserving inner nodes, finite leaves and infinite leaves.
\end{theorem}

Actually, the bijection between finite leaves of $\th$ and $\ttt$ is  a classical fact. To any $F\in\irr(\kh)$ we can associate a valuation $v_F\in\lfin(\th)$ defined as
$$
v_F(q)=\vb(q(\t)) \quad\mbox{ for all }q\in\khx,
$$
where $\t\in \kb$ is a root of $F$. By the henselian property, this construction does not depend on the choice of $\t$.
Clearly, $\supp(v_F)=F\khx$.

We denote the restriction of $v_F$ to $\kx$ by:
$$
w_F:=\left(v_F\right)_{\mid \kx}\in\lfin(\ttt).
$$
Now, $\supp(w_F)=N_{\kh/K}(F)\kx$, where  $N_{\kh/K}(F)\in\irr(K)$ is the monic generator of the prime ideal  $\left(F\khx\right)\cap\kx$.


\begin{proposition}\cite[Sec. 17]{endler}\label{lfinh}
The following two mappings are bijective:
$$
\irr(\kh)\lra\lfin(\th)\lra\lfin(\ttt),\qquad F\ \mapsto\ v_F\ \mapsto\ w_F.
$$
\end{proposition}

More generally, for any polynomial $g\in\irr(K)$, this construction facilitates the description of the extensions of $v$ to the simple extension $\kx/(g)$.
 
Since $K^h/K$ is a separable extension, we have $g=G_1\cdots G_r$ with pairwise different  $G_i\in \irr(\kh)$. Since $N_{\kh/K}(G_i)=g$ for all $i$,  
each $w_{G_i}\in\lfin(\ttt)$ induces  a valuation
$\wb_{G_i}$ on the field $K_g$.

\begin{theorem}\cite[Sec. 17]{endler}\label{endler}
	The extensions  of $v$ to $\kx/(g)$ are $\wb_{G_1},\dots,\wb_{G_r}$. 
\end{theorem}

In particular, $w_{G_1},\dots w_{G_r}$ are all finite leaves of $\ttt$ with support $g\kx$.

The aim of the OM-algorithm is to  compute, for each irreducible factor $G\in\khx$ of $g$, a chain of valuations in $\tinn$ getting ``sufficiently close" to the valuation $w_G$.

\section{Graded algebra and key polynomials}\label{secKP}
Take any $\mu\in\ttt$ and let $\p=\supp(\mu)$. For all $\alpha\in\g_\mu$, consider the $\oo$-modules:
$$
\ppa=\{g\in \kx\mid \mu(g)\ge \alpha\}\supset
\ppa^+=\{g\in \kx\mid \mu(g)> \alpha\}.
$$    
The \emph{graded algebras} of $v$ and $\mu$ are the integral domains:
$$
\gg_v=\bigoplus\nolimits_{\alpha\in\g}\left(\ppa\cap K\right)/\left(\ppa^+\cap K\right),\qquad \ggm=\bigoplus\nolimits_{\alpha\in\g_\mu}\ppa/\ppa^+.
$$
There is an obvious embedding of graded algebras $\gg_v\hk\ggm$.

Consider the \emph{initial term} mapping $\inm\colon \kx\to \ggm$, given by $\inm\p=0$ and 
$$
\inm g= g+\pset_{\mu(g)}^+, \quad\mbox{if }\ g\in \kx\setminus\p.
$$
We denote the \emph{grade} of $\inm g$ by $\op{gr}_\mu(\inm g)=\mu(g)$ as well.

Let us denote the set of all nonzero homogeneous elements in $\ggm$ by 
$$
\hmg=\left\{\inm g\mid g\in\kx\setminus\p\right\}.
$$
Let $\hmgu\subset\hmg$ be multiplicative group of all homogeneous units. \e

\defn Let $g,\,h\in \kx$.

We say that $g,h$ are \emph{$\mu$-equivalent}, and we write $g\smu h$, if $\inm g=\inm h$. 

We say that $g$ is \emph{$\mu$-divisible} by $h$, and we write $h\mmu g$, if $\inm h\mid \inm g$ in $\ggm$.

We say that $g$ is $\mu$-irreducible if $\inm g$ is a prime element. 

We say that $g$ is $\mu$-minimal if $g\nmid_\mu f$ for all nonzero $f\in \kx$ with $\deg(f)<\deg(g)$.\e

Recall that $\pi\in\hmg$ is a prime element if the homogeneous principal ideal of $\ggm$ generated by $\pi$ is a prime ideal. In this case, for all $t\in\hmg$ the order $n=\ord_\pi(t)$ is determined by the conditions  $\pi^n\mid t$, $\pi^{n+1}\nmid t$.

For all $\phi\in\kx\setminus K$ we define the \emph{truncation} $\mu_\phi$ as follows: 
$$
g=\sum\nolimits_{0\le n}a_n \phi^n,\quad \deg(a_n)<\deg(\phi)\quad\imp\quad \mu_\phi(g)=\min_{0\le n}\left\{\mu\left(a_n \phi^n\right)\right\}.
$$

This function $\mu_\phi$ is not necessarily a valuation, but it is useful to characterize the  $\mu$-minimality of $\phi$. Let us recall \cite[Prop. 2.3]{KP}.

\begin{lemma}\label{minimal0}
A polynomial  $\phi\in \kx\setminus K$ is $\mu$-minimal if and only if $\mu_\phi=\mu$.
\end{lemma}

\defn A  \emph{(Maclane-Vaqui\'e) key polynomial} for $\mu$ is a monic polynomial in $\kx$ which is simultaneously  $\mu$-minimal and $\mu$-irreducible. 
The set of key polynomials for $\mu$ is denoted $\kpm$. All  key polynomials are irreducible in $\kx$. 
\e

The existence of key polynomials characterizes the inner nodes of $\ttt$.

\begin{theorem}\cite[Thm. 4.4]{KP}\label{leaves}
	A valuation $\mu\in\ttt$ is a leaf if and only if  $\,\kpm=\emptyset$. This is equivalent to $\hmg=\hmgu$ too.
\end{theorem}

From now on, we assume that $\mu$ is a inner node of $\ttt$ and $\phi\in\kpm$ is a key polynomial. Also, we denote $$
\pi=\inm\phi,\qquad \ab=\inm a\ \mbox{ for all }a\in\kx.
$$ 

Since $\phi$ is $\mu$-minimal, 
Lemma \ref{minimal0} shows that $\mu$ acts on $\phi$-expansions as follows
\begin{equation}\label{muf}
g=\sum\nolimits_{0\le n}a_n \phi^n,\ \ \deg(a_n)<\deg(\phi)\ \imp\ \mu(g)=\min_{0\le n}\left\{\mu\left(a_n \phi^n\right)\right\}.
\end{equation}

Let $S_{\mu,\phi}(g):=\left\{0\le n\mid \mu\left(a_n \phi^n\right)=\mu(g)\right\}$. For all nonzero $g\in\kx$, we have
\begin{equation}\label{ordpi}
\overline g=\sum\nolimits_{n\in S_{\mu,\phi}(g)}\ab_n\pi^n,\qquad \ord_\pi(\overline g)=\min\left(S_{\mu,\phi}(g)\right).
\end{equation}
If $\phi$ is a key polynomial of minimal degree, then all these coefficients $\ab_n$ derived from a $\phi$-expansion are homogeneous units in $\ggm$ \cite[Prop. 3.5]{KP}.

Let $\ggm^0\subset\ggm$ be the subalgebra generated by the set of all homogeneous units. Equivalently, $\ggm^0$ is the relative algebraic closure  of $\gg_v$ in the embedding $\gg_v\hk\ggm$.

The following result is classically known (cf. for instance \cite[Prop. 4.5]{N2021}). 

\begin{theorem}\label{g0gm}
	Let $\phi$ be a key polynomial of minimal degree for $\mu$. 
	Then, the prime $\pi=\inm\phi$ is transcendental over $\ggm^0$ and
	$\ggm=\ggm^0[\pi]$.
\end{theorem}

\defn
For all nonzero $g\in\kx$ we define its $\mu$-\emph{degree} $\dgm(g)\in\N$ and \emph{leading coefficient} $\lcm(g)\in\ggm^0$ as the degree and leading coefficient of $\inm g$ as a polynomial in $\pi=\inm \phi$ with coefficients in $\ggm^0$, for some $\phi\in\kpm$ of minimal degree. \e

These definitions are independent of the choice of $\phi$ among all key poynomials of minimal degree for $\mu$.   
Note that an homogeneous element $\inm g$ is a unit if and only if $\dgm(g)=0$. Also, we have in general:  $\dgm(g)=\max\left(S_{\mu,\phi}(g)\right)$.\e

\defn
For an inner node $\mu\in\ttt$ we define its \emph{degree} as
$\deg(\mu)=\deg(\phi)$, 
where $\phi\in\kpm$ is a key polynomial of minimal degree. 

For a finite leaf $w_F\in\lfin$  we define $\deg(w_F)=\deg(F)$.

\subsection{Tangent directions of inner nodes}\label{subsecHPG}

A \emph{tangent direction} of an inner node $\mu$ of $\ttt$ is a $\mu$-equivalence class $\cl{\phi}\subset\kpm$ containing all key polynomials having the same initial term in $\ggm$. We denote the set of all tangent directions of $\mu$ by:
$$
\tdm=\kpm/\!\smu.
$$
This terminology is justified by item (ii) of the following result.

\begin{lemma}\cite[Lem. 2.2, Prop. 2.4]{VT}\label{propertiesTMN}
	Let $\mu<\nu$ be two nodes in $\ttt$. Let $\tmn$ be the set of monic polynomials $\phi\in\kx$ of minimal degree satisfying $\mu(\phi)<\nu(\phi)$.
	
	
	\,(i) \,The set $\tmn$ is a tangent direction of $\mu$.
	Moreover, for any $\phi\in\tmn$ and all nonzero $g\in\kx$, the equality $\mu(g)=\nu(g)$ holds if and only if $\phi\nmid_{\mu}g$.   
	
	
	(ii) \,If $\mu<\nu'$ for some $\nu'\in\ttt$, then
	
	\qquad\qquad\qquad$\tmn=\ty(\mu,\nu')\sii (\mu,\nu\,]\cap (\mu,\nu'\,]\ne\emptyset$. 
\end{lemma}

Let $\hpg\subset\hmg$ be the subset of all homogeneous prime elements in $\ggm$. 

The multiplicative group $\hmgu$  acts on $\hmg$ and $\hpg$ by ordinary multiplication. Let us denote the orbit of any $t\in\hmg$ by $$[\,t\,]_{\unit}=t\,\hmgu\in\hmg/\hmgu.$$ 


Clearly,  $\hpg/\hmgu$ can be identified with the set of all homogeneous principal prime ideals. 
The next result, which  follows easily from \cite[Thm. 6.8]{KP}, shows that all these prime ideals are generated by initial terms of key polynomials.

\begin{theorem}\label{uniquefact}
Let  $\mu\in\tinn$. 
\begin{enumerate}
\item [(i)] All $t\in\hmg$ factorize as a product of prime elements. The factorization is unique up to reordering the factors and multiplication by homogeneous units.
\item [(ii)] There is a canonical bijection:
$$
\tdm\lra \hpg/\hmgu,\qquad [\phi]_\mu\longmapsto [\inm\phi]_{\unit}.
$$
\end{enumerate}
\end{theorem}

For a given  $g\in\irr(K)$, let $\fff$ be the set of monic irreducible factors of $g$ in $\khx$. Also, for all $\mu\in\tinn$, let us denote 
$$\ffm=\left\{G\in\fff\mid \mu<w_G\right\}.$$

As mentioned in the last section, the OM-algorithm computes, for each $G\in\fff$, a chain of valuations in $\tinn$ getting sufficiently close to the valuation $w_G$.  

To this purpose, for a given valuation $\mu\in\tinn$, we need to compute the tangent directions of $\mu$ ``pointing out" to leaves $w_G\in\lfin$ associated to some $G\in\fff$; that is, the tangent directions of $\mu$ determined by the set $\ffm$.
These tangent directions are determined by the following criterion of Barnab\'e-Novacoski \cite[Thms. 1.1,1.3]{BN}.

\begin{theorem}\label{BN}
	Let $\mu\in\tinn$  and $g\in\irr(K)$. The image of the composition
	$$
	\ffm\to\tdm\to \hpg/\hmgu,\qquad G\mapsto \ty(\mu,w_G)\mapsto \ty(\mu,w_G)\hmgu.
	$$
is the set of prime homogeneous factors of $\,\inm g\in\ggm$, up to units.
	In particular, $\ffm=\emptyset$ if and only if $\,\inm g$ is a unit in $\ggm$.
\end{theorem}

We are led to the resolution of the following problem: 

{\bf \small  Given a nonzero $g\in\kx$,  compute the prime factorization of $\inm g$ in $\ggm$.}\e

After Theorem \ref{g0gm}, this amounts to factorize $\inm g$ in the algebra $\ggm^0[X]$ (where $X$ is an indeterminate). However, working in this algebra is computationally painful. A crucial feature of the OM-algorithm is that it provides the factorization of $\inm g$ by working in the subring $$\dm:=\pset_0/\pset_0^+\subset\ggm$$ of all homogeneous elements of grade zero, which is a polynomial ring with coefficients in a field.
This is the aim of Section \ref{subsecResPol}. 

\subsection{Residual polynomial operators}\label{subsecResPol}

Let  $\kappa=\kappa(\mu)$ be the relative algebraic closure of $k$ in $\km$.
There are canonical injective ring homomorphisms $$ k\hookrightarrow\kappa\hookrightarrow\dm\hookrightarrow \km.$$

\defn
Let $\gm^0=\left\{\mu(a)\mid a\in \kx,\ 0\le \deg(a)<\deg(\mu)\right\}$ be the subgroup of all grades of homogeneous units. By (\ref{muf}), we have $\gm=\gen{\gm^0,\mu(\phi)}$.

The \emph{relative ramification index} of $\mu$ is defined as $\,e=\erel(\mu)=\left(\gm\colon \gm^0\right)$.
Thus, $e$ is the least positive integer such that $e\mu(\phi)\in\gm^0$.\e

The following result is classical. A proof can be found in \cite[Thms. 4.5, 4.6]{KP}.


\begin{theorem}\label{Delta}
Let  $\pi=\inm\phi$ for a key polynomial $\phi$ of minimal degree. Take any homogeneous unit  $u\in\hmgu$ of grade $e\mu(\phi)$. Then, 
$\xi=\pi^eu^{-1}\in\dm$
is transcendental over $k$ and satisfies $\dm=\kappa[\xi]$.

Moreover, the canonical embedding $\dm\hookrightarrow \km$ induces an isomorphism $\ka(\xi)\simeq \km$. 
\end{theorem}

The pair $\phi, u$ determines a \emph{residual polynomial operator} 
$$
R=R_{\mu,\phi,u}\colon\;\kx\lra \kappa[y].
$$
Let us recall its definition. 
 We agree that $R(0)=0$.
For a nonzero $g \in K[x]$ with $\phi$-expansion $g=\sum\nolimits_{0\le n}a_n\phi^n$, let us denote
$$
S=S_{\mu,\phi}(g),\quad 
\ell_0=\min(S),\quad \ell=\max(S)=\dgm(g). 
$$
Note that $\bar{a}_\ell=\lcm(g)$. If we denote $\ga=\mu(\phi)$, then for all  $n\in\N$ we have
$$
n\in S\ \sii\  \mu(a_n)+n\ga=\mu(a_{\ell_0})+\ell_0\ga\ \sii\  (n-\ell_0)\ga=\mu(a_{\ell_0})-\mu(a_n).
$$
This implies that $(n-\ell_0)\ga$ belongs to $\gm^0$, so that  $n-\ell_0=je$ for some $j\in\N$. Since $\ell\in S$, this shows in particular that $\ell-\ell_0=de$ for some $d\in\N$. Let us denote
$$
\ell_j=\ell_0+je,\qquad  0\le j\le d.
$$
Note that $\ell_d=\ell$. Finally,  for all $0\le j\le d$, consider the \emph{residual coefficient}
\begin{equation}\label{resCoeff}
\zeta_j=\begin{cases}
\left(\bar{a}_\ell\right)^{-1}u^{j-d}\,\bar{a}_{\ell_j}\in \dm^*=\kappa^*,&\quad \mbox{ if }\ \ell_j\in S,\\
0,&\quad \mbox{ otherwise}.
\end{cases}
\end{equation}

\defn
$\,R(g)=\zeta_0+\zeta_1\,y+\cdots+\zeta_{d-1}y^{d-1}+y^d\in \kappa[y]$.\e

Since $\ell_0\in S$, we have $\zeta_0\ne 0$. Let us display
the essential property of this operator.

\begin{theorem}\label{Hmug}
For all nonzero $g\in\kx$,  $\,\inm g=\lcm(g)\, u^d\,\pi^{\ell_0}\,R(g)(\xi)$.
\end{theorem}

Indeed, if we denote $\ep=\lcm(g)=\bar{a}_\ell$, this follows immediately from:

$$
\ep^{-1}\inm g=\sum_{\ell_j\in S}\ep^{-1}\bar{a}_{\ell_j}\,\pi^{\ell_j}=\pi^{\ell_0}\sum_{\ell_j\in S}\ep^{-1}\bar{a}_{\ell_j}\,\pi^{je}=u^d\pi^{\ell_0}\sum_{\ell_j\in S}\zeta_j\,(\pi^e/u)^j.
$$



\begin{corollary}\cite[Cor. 5.4]{KP}\label{prodR} \ For all $g,h\in\kx$ we have $R(gh)=R(g)R(h)$.
\end{corollary}

With this tool in hand, \cite[Props. 6.3, 6.6]{KP} determine the whole set $\kpm$.

\begin{theorem}\label{charKP}
For a residually transcendental  $\mu$, take $\phi\in\kpm$ of minimal degree $m$. 
A monic $Q\in\kx$ is a key polynomial for $\mu$ if and  only if either
\begin{itemize}
\item $\deg(Q)=m$ \,and\; $Q\smu\phi$, or\vskip0.1cm
\item $\deg(Q)=me\deg(R(Q))$ \,and\; $R(Q)$ is irreducible in $\kappa[y]$.
\end{itemize}

Moreover, for all $Q,\,Q'\in\kpm$, we have
$$
Q\mmu Q'\sii Q\smu Q'\sii R(Q)=R(Q')\, \imp\, \deg(Q)=\deg(Q').
$$
\end{theorem}

\begin{corollary}\label{td} 
Let $\mu$ be a valuation on $\kx$ admitting a key polynomial $\phi\in\kpm$. For any valuation $\mu<\nu$, we have$$	\ty(\mu,\nu)=[\phi]_\mu\ \sii\ \mu(\phi)<\nu(\phi).	$$
\end{corollary}

\begin{proof}	
If $\ty(\mu,\nu)=[\phi]_\mu$, then $\mu(\phi)<\nu(\phi)$ by the definition of the tangent direction. 	Conversely, suppose $\mu(\phi)<\nu(\phi)$ and let $\ty(\mu,\nu)=[\varphi]_\mu$. Then, $\varphi\mmu \phi$, and this implies $\varphi\smu\phi$ by Theorem \ref{charKP}.  
\end{proof}\e

 It is easy to design a \emph{lifting routine}   \cite[Cor. 5.6]{KP}
$$
\op{lift}_{\mu,\phi}\colon \irr(\kappa)\setminus\{y\}\lra \kpm,\qquad \psi\longmapsto Q,
$$
to construct $Q\in\kpm$ with a prefixed  $R(Q)=\psi$.
We deduce a bijection  
$$
\tdm\lra \irr(\kappa),\qquad [Q]_\mu\longmapsto 
\begin{cases}
y,&\mbox{ if }Q\smu\phi,\\
R(Q),&\mbox{ otherwise},
\end{cases}
$$
which  depends on the choice the pair $\phi,u$. The variation of $R(Q)$ with respect to the pair $\phi,u$ is exhaustively discussed in \cite[Sec. 5]{KP}.

The factorization of $\inm g$ follows from Theorem \ref{Hmug} too. Let us factorize $R(g)$ as a product of powers of pairwise different irreducible polynomials in $\kappa[y]$:
$$
R(g)=\psi_1^{n_1}\dots\psi_r^{n_r},\qquad \psi_1,\dots,\psi_r\in\irr(\kappa).
$$
By \cite[Lem. 6.1]{KP}, we obtain a factorization $R(g)(\xi)=\psi_1(\xi)^{n_1}\dots\psi_r(\xi)^{n_r}$ as a product of homogeneous prime elements in $\ggm$.

Take $Q_i\in\kpm$ lifting $\psi_i$ and denote $\pi_i=\inm Q_i$, for all $i$. By Theorem \ref{Hmug},
$$
\pi_i\sim_{\unit} \psi_i(\xi)\quad\mbox{ for all }1\le i \le r,
$$
where $\sim_{\unit}$ indicates equality up to multiplication by some unit. 

Therefore, we obtain the following factorization of $\inm g$:
\begin{equation}\label{piFactors}
\inm g\sim_{\unit}\pi^{\ell_0}\psi_1(\xi)^{n_1}\dots\psi_r(\xi)^{n_r}\sim_{\unit}\pi^{\ell_0}\pi_1^{n_1}\dots\pi_r^{n_r}.
\end{equation}
The exponents $n_1,\dots,n_r$ are all positive, but $\ell_0=\min(S_{\mu,\phi}(g))$ might vanish.\e

\nn{\bf Convention. }Throughout the paper, we shall denote the operator $R_{\mu,\phi,u}$ simply by $R_{\mu,\phi}$, omitting its dependence on  the choice of a suitable homogeneous  unit $u\in\ggm$. Note that the degree of  $R_{\mu,\phi}(g)$ does not depend on the choice of $u$.

\section{\mlv chains}\label{secMLV}
\subsection{Depth-zero valuations and ordinary augmentations}\label{subsecDepth0}

For all $a\in K$, $\ga\in\gq\infty$, we may construct the \emph{depth-zero} valuation  $\mu=[v;\,x-a,\ga]\in\ttt$,
defined in terms of $(x-a)$-expansions as
$$
g=\sum\nolimits_{0\le n}a_n(x-a)^n\imp \mu(g)=\min\{v(a_n)+n\ga\mid 0\le n\}.
$$
Note that $\mu(x-a)=\ga$. 
If $\ga<\infty$, then $\mu$ is an inner node of $\ttt$ and $x-a$ is a key polynomial for $\mu$ of minimal degree. 
If $\ga=\infty$, then $\mu$ is the unique finite leaf of $\ttt$ with support $(x-a)\kx$. In both cases, $\deg(\mu)=1$.



Let $\mu\in\tinn$ be an inner node of $\ttt$.
For all $\phi\in\kpm$, $\ga\in\gq\infty$ such that $\mu(\phi)<\ga$, we may construct the \emph{ordinary} augmented valuation $\nu=[\mu;\,\phi,\ga]\in\ttt$,
defined in terms of $\phi$-expansions as
$$
g=\sum\nolimits_{0\le n}a_n\phi^n,\quad\deg(a_n)<\deg(\phi)\ \imp\ \nu(g)=\min\{\mu(a_n)+n\ga\mid 0\le n\},
$$
Note that $\nu(\phi)=\ga$, $\mu<\nu$ and  $\ty(\mu,\nu)=[\phi]_\mu$. 

If $\ga<\infty$, then  $\nu$ is an inner node of $\ttt$ and $\phi$ is a key polynomial for $\nu$ of minimal degree \cite[Cor. 7.3]{KP}. If $\ga=\infty$, then $\nu$ is a finite leaf of $\ttt$ with support $\phi\kx$. In both cases, $\deg(\nu)=\deg(\phi)$.


\subsection{Limit augmentation of valuations}\label{subsecLimAugm}
Let $\cc=\cfa$ be a totally ordered family of inner nodes of $\ttt$, not admitting a last element.
Assume that $A$ is a totally ordered set and $\ri<\rj$ if and only if $i<j$ in $A$.

We say that $\cc$ has \emph{stable degree}   if $\deg(\rho_{i})$ is stable for all sufficiently large $i\in A$. In this case, we denote this stable degree by $\dg(\cc)$.

We say that $g\in\kx$ is \emph{$\cc$-stable} if for some index $i\in A$, we have
$\ri(g)=\rj(g)$ for all $j>i$.
We may define a \emph{stability function} $\rhc(g)=\max\{\ri(g)\mid i\in A\}$,
on the set of all $\cc$-stable polynomials.\e

\defn
A \emph{limit key polynomial} for $\cc$ is a monic $\cc$-unstable polynomial of minimal degree. Let $\kpi(\cc)$ be the set of all limit key polynomials. Since the product of stable polynomials is stable, all limit key polynomials are irreducible in $ \kx$.

We say that $\cc$ is an \emph{essential continuous family} of valuations if it has stable degree and admits limit key polynomials of degree  greater than $\dg(\cc)$.\e

For all $\phi\in\kpi\left(\cc\right)$, $\ga\in\gq\infty$ such that $\ri(\phi)<\ga$ for all $i\in A$, we may construct the \emph{limit augmented} valuation $\mu=[\cc;\,\phi,\ga]\in\ttt$,
defined in terms of $\phi$-expansions as: 
$$ g=\sum\nolimits_{0\le n}a_n\phi^n,\quad\deg(a_n)<\deg(\phi)\ \imp\ \mu(g)=\min\{\rhc(a_n)+n\ga\mid 0\le n\}.
$$
Note that $\mu(\phi)=\ga$ and $\ri<\mu$ for all $i\in A$.

If $\ga<\infty$, then $\mu$ is an inner node of $\ttt$ and $\phi$ is a key polynomial for $\mu$ of minimal degree \cite[Cor. 7.13]{KP}.
If $\ga=\infty$, then $\mu$ is a finite leaf of $\ttt$ with support $\phi\kx$. In both cases, $\deg(\mu)=\deg(\phi)$.

\subsection{Maclane--Vaqui\'e chains}

Take a chain of finite length $r$, of valuations in $\ttt$
\begin{equation}\label{depthMLV}
v\ \stackrel{\phi_0,\ga_0}\lra\  \mu_0\ \stackrel{\phi_1,\ga_1}\lra\  \mu_1\ \stackrel{\phi_2,\ga_2}\lra\ \cdots
\ \lra\ \mu_{r-1} 
\ \stackrel{\phi_{r},\ga_{r}}\lra\ \mu_{r}=\mu
\end{equation}
in which $\mu_0=[v;\,\phi_0,\ga_0]$ is a depth-zero valuation, and each other node is an augmentation  of the previous node, of one of the two types:\e

\emph{Ordinary augmentation}: \ $\mu_{n+1}=[\mu_n;\, \phi_{n+1},\ga_{n+1}]$, for some $\phi_{n+1}\in\kp(\mu_n)$.\e

\emph{Limit augmentation}:  \ $\mu_{n+1}=[\cc_n;\, \phi_{n+1},\ga_{n+1}]$,  for some $\phi_{n+1}\in\kpi(\cc_n)$, where $\cc_n$ is an essential continuous family  whose first valuation is $\mu_n$.\e

\defn
A chain of mixed augmentations as in (\ref{depthMLV}) is said to be  a \emph{\mlv (MLV) chain}  if every augmentation step satisfies:
\begin{itemize}
	\item If $\,\mu_n\to\mu_{n+1}\,$ is ordinary, then $\ \deg(\mu_n)<\deg(\mu_{n+1})$.
	\item If $\,\mu_n\to\mu_{n+1}\,$ is limit, then $\ \deg(\mu_n)=\deg(\cc_n)$ and $\ \phi_n\not \in\ty(\mu_n,\mu_{n+1})$. 
\end{itemize}\e

In this case, we have $\mu(\phi_n)=\ga_n$ for all $n$. 
As shown in \cite[Sec. 4.1]{MLV}, the MLV chain induces a chain of value groups 
$$
\g_{\mu_{-1}}:=\g\subset\g_{\mu_0}\subset \cdots \subset \g_{\mu_{r}}= \gm,
$$
such that $\g_{\mu_{n-1}}=\g_{\mu_n}^0$ for all $0\le n\le r$, and 
\begin{equation}\label{valuegroups}
\g_{\mu_{n}}=\gen{\g_{\mu_{n-1}},\ga_{n}} \ \mbox{ if }\ \ga_n<\infty,\qquad 
\gm=\g_{\mu_{r-1}} \ \mbox{ if }\ \ga_r=\infty.
\end{equation}
As shown in \cite[Sec. 5.1]{MLV}, the MLV chain induces a tower of finite and simple extensions of fields
$$
\ka(\mu_{-1}):=k\,\to\,\ka(\mu_0)\,\to\,\cdots\,\to\,\ka(\mu_r)=\ka(\mu).
$$
For all $0\le n\le r$, let us denote
\begin{equation}\label{esfs}
e_n=\left(\g_{\mu_n}\colon\g_{\mu_{n-1}}\right)=\erel(\mu_n),\quad f_n=\left[\ka(\mu_n)\colon\ka(\mu_{n-1})\right]=\deg\left(R_{\mu_{n-1},\phi_{n-1}}(\phi_n)\right),
\end{equation}
the last equality by \cite[Lem. 5.2,5.3]{MLV}. 

If $\mu$ has nontrivial support $g\kx$, then we can read in the MLV chain of $\mu$ the ramification index $e(\overline{\mu}/v)$ and residual degree $f(\overline{\mu}/v)$ of the valuation $\overline{\mu}$ induced by $\mu$ on the field $\kx/(g)$. Obviously, $\g_{\overline{\mu}}=\gm$ and (by definition) $k_{\overline{\mu}}=\km$.

\begin{proposition}\cite[Thm. 5.4]{MLV}\label{infinity}
If $\ga_r=\infty$, then $\mu$ is a finite leaf with $k_\mu=\kam$. In particular, $e(\overline{\mu}/v)=e_0\cdots e_{r-1}$ and $f(\overline{\mu}/v)=f_0\cdots f_{r}$.
\end{proposition}

The following theorem is due to Maclane, for the discrete rank-one case  \cite{mcla}, and Vaqui\'e for the general case  \cite{Vaq}. Another proof may be found in  \cite[Thm. 4.3]{MLV}. 

\begin{theorem}\label{main}
All $\mu\in\tinn\sqcup\lfin$ are the end node of a finite MLV chain.
\end{theorem}


The main advantage of MLV chains is that they are essentially unique, so that we may read in them several data intrinsically associated to the valuation $\mu$.
For instance, the sequence $\left(\deg(\mu_n)\right)_{n\ge0}$, the character ``ordinary" or ``limit" of each augmentation step, and the numerical data from (\ref{esfs}) are intrinsic features of $\mu$ \cite[Sec. 4.3]{MLV}.  \e

\defn
The \emph{depth} of $\mu$ is the length of any MLV chain with end node  $\mu$.

We say that $\mu$ is \emph{inductive} if all augmentations in its MLV chain are ordinary. 

\subsection{Inductive valuations and henselization}\label{subsecIndRig}
Let $(\kh,\vh)$ be a henselization of $(K,v)$ and let $\muh$ be the unique common extension of $\mu$ and $\vh$ to $\khx$  (Theorem \ref{Rig}). The mapping $\inm g\ \mapsto\ \inmh g$, for all $g\in\kx$, induces an embedding $\ggm\hk\gg_{\muh}$ of graded algebras.

\begin{theorem}\cite{NN}\label{NN}
For all	valuations $\mu$ on $\kx$ the canonical embedding $\ggm\hk\ggmh$ is an isomorphism of graded algebras.
\end{theorem}

\begin{lemma}\label{IndRig}
	If $\mu$ be an inductive valuation on $\kx$, admitting a MLV chain as in (\ref{depthMLV}). Then, $\kpm\subset\kp(\muh)$ and  $\muh$ is inductive, admitting a MLV chain with the same length $r$ and data $(\phi_n,\ga_n)$ 	for all $\,0\le n\le r$:
	$$
	\vh\ \stackrel{\phi_0,\ga_0}\lra\  \muh_0\ \stackrel{\phi_1,\ga_1}\lra\  \muh_1\ \stackrel{\phi_2,\ga_2}\lra\ \cdots
	\ \lra\ \muh_{r-1} 
	\ \stackrel{\phi_{r},\ga_{r}}\lra\ \muh_{r}=\muh.
	$$
Moreover, the numerical data $e_0,\dots,e_r;\, f_0,\dots,f_r$  attached to both chains coincide. 
\end{lemma}

\begin{proof}
The first statement is well-known. It can be found, for instance, in \cite{NMO}.

The value groups of both MLV chains coincide by (\ref{valuegroups}). Hence, both chains determine the same data   $e_0,\dots,e_r$. Finally, the two towers of fields $\ka(\mu_n)\to\ka(\mu_{n+1})$ and $\ka(\mu_n^h)\to\ka(\mu_{n+1}^h)$ are isomorphic by Theorem \ref{NN}. Hence, both chains determine the same data   $f_0,\dots,f_r$.
\end{proof}

\section{Newton polygons}\label{secNewton}
Consider two points $P=(n,\alpha),\ Q=(m,\beta)$ in the $\Q$-vector space $\gqq$.
The segment joining $P$ and $Q$ is the subset
$$
S=\left\{P+ \dta\,\overrightarrow{PQ}\mid\  \dta \in  \Q,\ 0\le \dta\le 1\right\}\subset\gqq.
$$
If $n\ne m$, this segment has a natural slope
$$
\slp(S)=(\beta-\alpha)/(m-n)\in\gq.
$$


A subset of $\gqq$ is \emph{convex} if it contains the segment joining any two points in the subset.
The \emph{convex hull} of a finite subset $C\subset \gqq$ is the smallest convex subset of $\gqq$ containing $C$. 

The border of this hull is a sequence of chained segments.  If the points in $C$ have different abscissas, the leftmost and rightmost points are joined by two different chains of segments along the border, called the \emph{upper} and \emph{lower} convex hull of $C$.

\subsection{Classical Newton polygons}\label{subsecClassNP}
Let $\vb$ be a fixed extension of $v$ to $\kb$ and $(\kh,\vh)$ the corresponding henselization of $(K,v)$ .

Let us recall the classical \emph{Newton polygon} operator
$$
N_{v,x}\colon\, \kx\lra \pset\left({\gqq}\right),
$$
where $\pset\left({\gqq}\right)$ is the power set of the rational vector space $\gqq$. 

The Newton polygon of the zero polynomial is the empty set. \e

\defn
For a nonzero $g=a_0+\cdots +a_\ell x^\ell\in \kx$, the Newton polygon $N_{v,x}(g)$ is the lower convex hull of the finite cloud of points 
$
\left\{\left(n,v(a_n)\right)\mid n\ge0\right\}
$.\e

Thus, $N:=N_{v,x}(g)$ is either a single point or a chain of segments, $S_1,\dots, S_r$, called the \emph{sides} of the polygon, ordered from left to right by increasing slopes.  

The abscissa of the left endpoint of $N$ is $\ord_x(g)$.


We define the \emph{length} $\ell(S_i)$ of a side as the length of its projection to the $x$-axis.

For all $g\in\kx$, let $Z(g)$ be the multiset of all roots of $g$ in $\kb$, counting multipli\-ci\-ties.
Also, let $V(g)$ be the multiset of all values $\vb(\t)\in\gq\infty$, for $\t$ running on $Z(g)$.
Both multisets have cardinality $\ell=\deg(g)$.

\begin{theorem}\label{classicalNP}
For a nonzero $g\in \kx$, suppose that the sides of $N_{v,x}(g)$ have slopes $-\la_1<\cdots<-\la_r$ and lengths $\ell_1,\dots,\ell_r$. Then,  $V(g)=\left\{\la_1^{\left(\ell_1\right)},\dots,\la_r^{\left(\ell_r\right)}\right\}$.    
\end{theorem}

By the henselian property, all roots of an irreducible polynomial in $\khx$ have the same $\vb$-value. Hence, Theorem \ref{classicalNP} determines a factorization:
$$
g=G_0G_1\cdots G_r,\quad G_i\in\khx,\quad \deg(G_i)=\ell_i,
$$
where $G_0=x^{\ord_x(g)}$ and, for $i\ge1$, $G_i$ is the product of all irreducible factors of $g$ in $\khx$ such that the $\vb$-value of its roots is $\la_i$. 

Therefore, the Newton polygon determines a \emph{dissection} of the multiset of all irreducible factors of $g$ in $\khx$, counting multiplicities.

\begin{figure}
\caption{Newton polygon $N=\np(g)$ of $g\in \kx$. }\label{figNmodel}
\begin{center}
\setlength{\unitlength}{4mm}
\begin{picture}(20,10)
\put(14.8,6){$\bullet$}\put(13.5,7){$\bullet$}\put(12.8,2){$\bullet$}\put(10,2){$\bullet$}\put(12,0.4){$\bullet$}\put(7.4,1.5){$\bullet$}
\put(6.2,8.35){$\bullet$}\put(3.75,6.9){$\bullet$}\put(2.8,8.35){$\bullet$}
\put(-1,3.6){\line(1,0){20}}\put(0,0){\line(0,1){10}}
\put(7.6,1.8){\line(-2,3){4.5}}\put(7.6,1.83){\line(-2,3){4.5}}
\put(7.6,1.8){\line(4,-1){4.7}}\put(7.6,1.83){\line(4,-1){4.7}}
\put(12.2,0.5){\line(1,2){2.8}}\put(12.2,0.53){\line(1,2){2.8}}
\multiput(3,3.5)(0,.25){21}{\vrule height2pt}
\multiput(15,3.5)(0,.25){11}{\vrule height2pt}
\put(1.7,2.6){\begin{footnotesize}$\ord_{\phi}(g)$\end{footnotesize}}
\put(14.5,2.6){\begin{footnotesize}$\ell(N)$\end{footnotesize}}
\put(18.5,2.6){\begin{footnotesize}$\Q$\end{footnotesize}}
\put(-1.4,9.2){\begin{footnotesize}$\gq$\end{footnotesize}}
\put(-.6,2.8){\begin{footnotesize}$0$\end{footnotesize}}
\end{picture}
\end{center}
\end{figure}
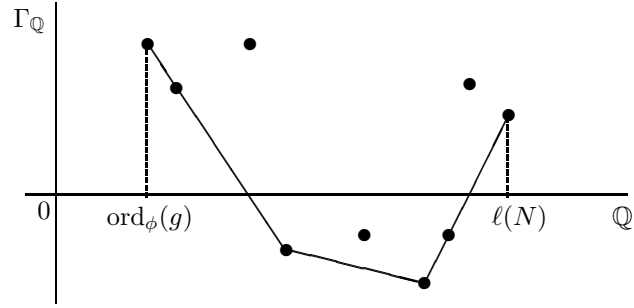

\subsection{General Newton polygons}\label{subsecGeneralNP}
A \emph{type}  is a pair $(\mu,\phi)$, where $\mu$ is an inner node of $\ttt$ and $\phi$ is a key polynomial for $\mu$.

Any  type $(\mu,\phi)$ yields a Newton polygon operator
$$
\np\colon\, \kx\lra \pset\left({\gqq}\right).
$$

The Newton polygon of the zero polynomial is the empty set. 
For a nonzero $g\in \kx$ with $\phi$-expansion 
$$g=\sum\nolimits_{0\le n}a_n\phi^n,\quad a_n\in\kx,\quad \deg(a_n)<\deg(\phi), $$
we define $N:=\np(g)$ as the lower convex hull of the finite set 
$
\left\{\left(n,\mu\left(a_n\right)\right)\mid n\ge0\right\}
$.

The abscissa of the left endpoint of $N$ is $\ord_{\phi}(g)$ in $\kx$.  The abscissa $$\ell(N)=\left\lfloor \deg(g)/\deg(\phi)\right\rfloor$$ of the right endpoint of $N$ is called the \emph{length} of $N$.
In Figure \ref{figNmodel} we display the typical shape of such a polygon.\e

\defn
For all $\la\in\gq$, the \mbox{\emph{$\la$-component}} $S_\la(N)\subset N$ is the intersection of $N$ with the line of slope $-\la$ which first touches $N$ from below. In other words,
$$S_\la(N)= \{(n,\al)\in N\,\mid\, \al+n\la\mbox{ is minimal}\,\}.$$
The abscissas of the endpoints of $S_\la(N)$ are denoted \ $n_{\la}\le n'_{\la}$.\e

If $N$ has a side $S$ of slope $-\la$, then $S_\la(N)=S$. Otherwise, $S_\la(N)$ is a vertex of $N$.  Figure \ref{figComponent0} illustrates both possibilities.\e

\defn
 $N$ is \emph{one-sided} of slope $-\la$, if 
$N=S_\la(N)$, $n_{\la}=0$ and $n'_{\la}>0$.

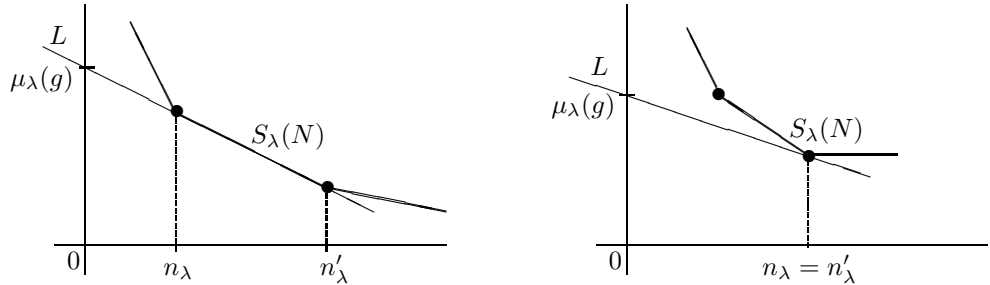
\begin{figure}
	\caption{$\la$-component of $N=\np(g)$. The line $L$ has slope $-\la$ and cuts the vertical axis at $(0,\mul(g))$, if $\la>\mu(\phi)$ and $\mul=[\mu;\,\phi,\la]$.
	}\label{figComponent0}
	\begin{center}
		\setlength{\unitlength}{4mm}
		\begin{picture}(30,9)
		\put(2.8,4.8){$\bullet$}\put(7.8,2.25){$\bullet$}
		\put(-0.2,6.5){\line(1,0){0.5}}
		\put(-1,0.6){\line(1,0){13}}\put(0,-0.4){\line(0,1){9}}
		\put(-1.4,7.2){\line(2,-1){11}}
		\put(3,5){\line(-1,2){1.5}}\put(3,5.04){\line(-1,2){1.5}}
		\put(3,5){\line(2,-1){5}}\put(3,5.04){\line(2,-1){5}}
		\put(8,2.5){\line(5,-1){4}}\put(8,2.54){\line(5,-1){4}}
		\multiput(3,.4)(0,.25){18}{\vrule height2pt}
		\multiput(8,.4)(0,.25){8}{\vrule height2pt}
		\put(7.8,-.4){\begin{footnotesize}$n'_{\la}$\end{footnotesize}}
		\put(2.6,-.4){\begin{footnotesize}$n_{\la}$\end{footnotesize}}
		\put(-1.2,7.3){\begin{footnotesize}$L$\end{footnotesize}}
		\put(-.6,-.2){\begin{footnotesize}$0$\end{footnotesize}}
		\put(-2.5,5.9){\begin{footnotesize}$\mul(g)$\end{footnotesize}}
		\put(5.5,4){\begin{footnotesize}$S_\la(N)$\end{footnotesize}}
		\put(20.8,5.35){$\bullet$}\put(23.8,3.3){$\bullet$}
		\put(17.75,5.58){\line(1,0){0.5}}
		\put(17,0.6){\line(1,0){13}}\put(18,-0.4){\line(0,1){9}}
		\put(16.1,6.2){\line(3,-1){10}}
		\put(24,3.6){\line(-3,2){3}}\put(24,3.64){\line(-3,2){3}}
		\put(21,5.8){\line(-1,2){1}}\put(21,5.84){\line(-1,2){1}}
		\put(24,3.6){\line(1,0){3}}\put(24,3.63){\line(1,0){3}}
		\multiput(24,.5)(0,.25){12}{\vrule height2pt}
		\put(22.5,-.4){\begin{footnotesize}$n_{\la}=n'_{\la}$\end{footnotesize}}
		\put(16.8,6.25){\begin{footnotesize}$L$\end{footnotesize}}
		\put(17.4,-.2){\begin{footnotesize}$0$\end{footnotesize}}
		\put(23.4,4.2){\begin{footnotesize}$S_\la(N)$\end{footnotesize}}
		\put(15.5,5){\begin{footnotesize}$\mul(g)$\end{footnotesize}}
		\end{picture}
	\end{center}
\end{figure}

\subsection{Dissection by Newton polygons}\label{subsecDissNP}\mbox{\null}\e

\defn
The \emph{principal Newton polygon} $\npp(g)$ is the polygon formed by the sides of $\np(g)$ of slope less than $-\mu(\phi)$. 

If $\np(g)$ has no sides of slope less than $-\mu(\phi)$, then $\np^+(g)$ is defined to be the left endpoint of $\np(g)$. 
\e

Clearly, $S_{\mu,\phi}(g)=\left\{0\le n\mid \mu\left(a_n \phi^n\right)=\mu(g)\right\}$ coincides with the set of abscissas of the points lying on the segment $S_{\mu(\phi)}(g)$. In particular, 
$\ell\left(\np^+(g)\right)=\min(S_{\mu,\phi}(g))$.
Hence, the following result is an immediate consequence of (\ref{ordpi}).

\begin{lemma}\label{length}
The integer $\ell\left(\npp(g)\right)$  is the order with which the prime element $\inm\phi$ divides $\inm g$ in the graded algebra $\ggm$.
\end{lemma}

\begin{lemma}\label{ordinate}
Let $(\mu,\phi)$ be a type. For $\la>\mu(\phi)$, let  $\mul=[\mu;\,\phi,\la]$.  Then, for all nonzero $g\in\kx$, the line of slope $-\la$ which first touches $\np(g)$ from below, cuts the vertical axis at the point of ordinate $\mul(g)$. 
\end{lemma}

\begin{proof}
This line cuts the vertical axis at the point with ordinate the common value of $\al+n\la$, for all $(n,\al)\in S_\la(\np(g))$ (cf. Figure \ref{figComponent0}).
\end{proof}\e

\defn
For a nonzero $g\in\kx$ and $\la\in\gq$, consider the multisets
$$
\ffmph=\{G\in\irr(\kh)\mid\ G\mid g,\ \mu< w_G\ \mbox{and }\ty(\mu,w_G)=[\phi]_\mu\}.
$$
$$
\ffmph(\la)=\{G\in\ffmph\mid\ w_G(\phi)=\la\}.
$$

Theorem \ref{conjNewton} below shows that  $\ffmph(\la)\ne\emptyset$ if and only if $-\la$ is one of the slopes of $\npp(g)$. As a consequence, we obtain a dissection
$$
\ffmph=\bigsqcup\nolimits_\la \ffmph(\la),
$$
for $-\la$ running on the slopes of $\npp(g)$.

Moreover, we get some information about the degrees and the number of irreducible factors of $g$ in these sets. 
The crucial point is the consideration of a special irreducible factor of $\phi$ in $\khx$, determined by the valuation $\mu$.\e

\defn
The valuation $[\mu;\,\phi,\infty]$ has support $\phi\kx$. As we saw in Section \ref{secComm}, there exists a unique irreducible factor $Q=Q_{\mu,\phi}\in\irr(\kh)$ of $\phi$ such that
$$
[\mu;\,\phi,\infty]=w_{Q}.
$$  
We say that $Q$ is the \emph{irreducible factor of $\phi$ over $\khx$ determined by} $\mu$.
	
\begin{theorem}\label{conjNewton}
Let $\mu$ be an inner node of $\ttt$ and $\phi\in\kpm$. Denote $N=\npp(g)$ and let $Q\in\irr(\kh)$ be the irreducible factor of $\phi$ determined by $\mu$. Then,
\begin{enumerate}
\item[(i)] All $G\in \ffmph$ have degree a multiple of $\deg(Q)$.
\item[(ii)] For all $\la\in\gq$, we have 
$$\sum\nolimits_{G\in\ffmph(\la)}\deg(G)=\ell\left(S_\la(N)\right) \deg(Q).$$
\end{enumerate}
In particular, if $\ell\left(S_\la(N)\right)=1$, then $\ffmph(\la)$ contains a unique irreducible factor of $g$, and this factor has degree $\deg(Q)$. 
\end{theorem}

If $v$ has rank one, then $\phi=Q$ and this theorem follows easily from Montes' original arguments in the discrete rank-one case. The proof in the general case is much more involved. We postpone it to Section \ref{subsecNPH}, which is entirely devoted to this purpose. 

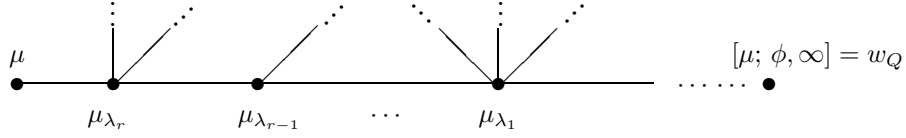
\begin{figure}
	\caption{Double dissection of $\ffmph$. The interval $[\mu,w_Q]$ consists of all augmentations $[\mu;\,\phi,\ga]$ for $\ga\in\gq\infty$ satisfying $\ga>\mu(\phi)$.}\label{figDD}
	\begin{center}
		\setlength{\unitlength}{4mm}
		\begin{picture}(25,5)
		\put(-2,1){$\bullet$}\put(-1.6,1.3){\line(1,0){21}}\put(23,1){$\bullet$}
		\put(6,1){$\bullet$}\put(1.2,1){$\bullet$}\put(14,1){$\bullet$}
		\put(-2,2){\begin{footnotesize}$\mu$\end{footnotesize}}
		\put(21.5,1){\begin{footnotesize}$\cdots$\end{footnotesize}}
		\put(20.2,1){\begin{footnotesize}$\cdots$\end{footnotesize}}
		\put(10,0){\begin{footnotesize}$\cdots$\end{footnotesize}}
		\put(22,2){\begin{footnotesize}$[\mu;\,\phi,\infty]=w_Q$\end{footnotesize}}
		\put(0.6,0){\begin{footnotesize}$\mu_{\la_{r}}$\end{footnotesize}}
		\put(5.6,0){\begin{footnotesize}$\mu_{\la_{r-1}}$\end{footnotesize}}
		\put(13.6,0){\begin{footnotesize}$\mu_{\la_{1}}$\end{footnotesize}}
		\put(1.45,1.1){\line(0,1){2}}\put(1.35,1.1){\line(1,1){2}}
		\put(1.3,3.2){\vdots}\multiput(3.4,3)(0.25,.25){3}{$\cdot$}
		\put(6.25,1.2){\line(1,1){2}}\multiput(8.4,3.05)(0.25,.25){3}{$\cdot$}
		\put(14.25,1.1){\line(0,1){2}}\put(14.2,1.1){\line(1,1){2}}
		\put(14.2,3.2){\vdots}\multiput(16.4,3.1)(0.25,.25){3}{$\cdot$}
		\multiput(11.8,3.1)(-0.25,.25){3}{$\cdot$}\put(14.3,1.1){\line(-1,1){2}}
		\end{picture}
	\end{center}
\end{figure}

\subsection{Dissection by factorization of residual polynomials}\label{subsecDissRPol}

For $g\in\kx$ and  our fixed type $(\mu,\phi)$ as above, let $-\la_1<\cdots<-\la_r$ be the slopes of $\npp(g)$. 

Let us assume that $\phi\nmid g$ in $\kx$, so that $n_{\la_1}=0$.

For each slope $-\la$ of $\npp(g)$, consider the augmentation 
$$
\mul=[\mu;\,\phi,\la].
$$
The multiset $\ffmph(\la)$ can be further dissected by factorizing the residual polynomial $R_{\mul,\phi}(g)$  in $\kappa(\mul)[y]$: 
\begin{equation}\label{factR}
R_{\mul,\phi}(g)=\psi_1^{n_1}\cdots\psi_s^{n_s},\qquad\psi_1,\dots,\psi_s\in\irr(\kappa(\mul)).
\end{equation}

Consider arbitrary lifts $\varphi_i=\op{lift}_{\mul,\phi}(\psi_i)\in\kp(\mul)$, for all $1\le i\le s$.

As we saw in Section \ref{subsecResPol}, this leads to a factorization of $\inn_{\mul}g$ into a product of pairwise different homogeneous prime elements in $\gg_{\mul}$, up to units:
$$
\inn_{\mul}g\,\sim_{\unit}\,\pi^{n_\la}\pi_1^{n_1}\cdots\pi_s^{n_s}, \qquad \pi=\inn_{\mul}\phi, \quad \pi_i=\inn_{\mul}\varphi_i,\ 1\le i\le s. 
$$
By Theorem \ref{BN}, if $\la<\la_1$ (so that $n_\la>0$), the irreducible factors of $g$  in the set $\ff_{\mul}(g)$ determine exactly $s+1$ tangent directions of $\mul$, which are precisely
$$
[\phi]_{\mul},\ \, [\varphi_1]_{\mul},\dots,[\varphi_s]_{\mul}. 
$$
Take any irreducible factor $G\in\ff_{\mul}(g)$. By  Lemma \ref{propertiesTMN}, $w_G(\phi)>\la$ if and only if $\ty(\mul,w_G)=[\phi]_{\mul}$. 
Therefore, the irreducible factors $G$ with $w_G(\phi)=\la$ are distributed among the rest of tangent directions. We get a dissection
$$
\ffmph(\la)=\bigsqcup\nolimits_{i=1}^s \ff_{\mul,\varphi_i}(g).
$$
Finally, if $\la=\la_1$, then $n_\la=0$ and we get directly the same dissection.
Figure \ref{figDD} illustrates this ``double-dissection" process.

The computation of each principal Newton polygons $N^+_{\mul,\varphi_i}(g)$ leads to further double-dissections. By Lemma \ref{length}, we know a priori that
\begin{equation}\label{priori}
\ell\left(N^+_{\mul,\varphi_i}(g)\right)=n_i,\quad 1\le i\le s.
\end{equation}
This is relevant from a computational perspective: \emph{in all required computations of principal Newton polygons, we know a priori the length of the polygon}. Thus, we need only to implement the following truncated-expansions subroutine.\e

\nn{\bf Newton polygon}  {\sl NP}$(\mu,\phi,\ell)(g)$

\nn{\bf Input: }$g\in\kx$, \ $\mu\in\tinn$, $\phi\in \kpm$, $\ell\in\N$

\nn{\bf Output: }A list of all sides of $\npp(g)$\e

\nn compute the first $\ell+1$ coefficients $a_0,\dots,a_\ell$ of the $\phi$-expansion of $g$

\nn{\bf return} \;lower convex hull of the set $\{(n,\mu(a_n))\mid 0\le n\le \ell\}$. \e

\defn
We say that the type $(\mu,\phi)$ \emph{singles out} an irreducible factor of $g$ in $\khx$, if $\ffmph=\{G\}$ and $\deg(G)=\deg(Q_{\mu,\phi})$.\e

From (\ref{priori}) and Theorem \ref{conjNewton}, we derive another relevant observation.

\begin {proposition}\label{N=1}
If in the factorization (\ref{factR}) we have $n_i=1$, then  the pair $\left(\mul,\varphi_i\right)$ singles out an irreducible factor of $g$ in $\khx$. 
\end {proposition}

\subsection{Newton polygons and henselization}\label{subsecNPH}

\subsubsection{Addition of Newton polygons}\label{subsubsecAdd}

There is an addition law for Newton polygons. 
Consider two  polygons $N$, $N'$ with sides $S_1,\dots,S_r$, $S'_1,\dots,S'_{s}$, respectively. 

The left endpoint of the sum $N+N'$ is  the vector sum in $\Q\times\gq$ of the left endpoints of $N$ and $N'$, whereas the sides of $N+N'$ are obtained by joining to this endpoint all sides in the multiset $\left\{S_1,\dots,S_r,S'_1,\dots,S'_{s}\right\}$, ordered by increasing slopes.

If one of the polygons is a one-point polygon (say  $N'=\{P\}$), then it has an empty set of sides and the sum $N+N'$ coincides with the vector sum $N+P$ in $\Q\times \gq$.


\begin{theorem}\label{product}
	For all $\phi\in\kpm$ and nonzero $g,h\in \kx$, we have   
	$$
	\npp(gh)=\npp(g)+\npp(h).
	$$
\end{theorem}

This result is proved in \cite[Thm. 4.8]{RPO} for $\mu$ an inductive valuation. However, the proof is valid in the general case.

\subsubsection{Newton polygons with respect to henselian valuations}
We assume in this section that the valued field $(K,v)$ is henselian. 

The following result is crucial for our purpose. It was proved in \cite[Sec. 4]{defless} for inductive valuations. However, the proof is valid in the general case. 

\begin{theorem}\label{fundamental}
	Let $Q\in\kpn$ for some valuation $\nu$ on $\kx$. Then,
$$
	Q\mnu F\ \sii\ \nu<v_F \ \mbox{ and }\ \ty(\nu,v_F)=[Q]_\nu,
$$
for all $F\in\irr(K)$.
	Moreover, if these conditions hold, then:
	\begin{enumerate}
		\item[(i)] Either $F=Q$, or the Newton polygon $N_{\nu,Q}(F)$ is one-sided of slope $-v_F(Q)$.    
		\item[(ii)] $F\snu Q^\ell$ with $\ell=\ell(N_{\nu,Q}(F))=\deg(F)/\deg(Q)$. 
	\end{enumerate}
\end{theorem}



Let us rewrite Theorem \ref{conjNewton} in the henselian case and show that it follows immediately from Theorem \ref{fundamental}. 

\begin{theorem}\label{ConjHensel}
	For a henselian $(K,v)$, let $\nu$ be a valuation on $\kx$, $Q\in\kpn$,  $g\in\kx$  monic and $N=N^+_{\nu,Q}(g)$. Then,
	\begin{enumerate}
		\item [(i)] For all $G\in\ffnq$, $\deg(G)$ is a multiple of $\deg(Q)$.
		\item [(ii)] For all $\ep\in\gq$, we have
		$$\sum\nolimits_{G\in\ffnq(\ep)}\deg(G)=\ell\left(S_\ep(N)\right) \deg(Q).$$
	\end{enumerate}
\end{theorem}

\begin{proof}
	Let $\mathcal{F}(g)$ be the multiset of irreducible factors of $g$. By Theorem \ref{fundamental},
	$$
	\ffnq=\{G\in \mathcal{F}(g)\mid \ \nu<v_G,\ \ty(\nu,v_G)=[Q]_\nu\}=\{G\in\mathcal{F}(g)\mid \ Q\mnu G\},
	$$
	and all polynomials in this set have degree a multiple of $\deg(Q)$. This proves (i).
	
	For all $G\in\ffnq$, Theorem \ref{fundamental} shows that
	the Newton polygon $N_{\nu,Q}(G)$ is one-sided of length $\deg(G)/\deg(Q)$ and slope $-v_G(Q)$, with $v_G(Q)>\nu(Q)$. In particular,
	$$
	\ell\left(S_{v_G(Q)}(N^+_{\nu,Q}(G))\right)=\ell(N_{\nu,Q}(G))=\deg(G)/\deg(Q).
	$$
	
	Now, recall that $\ffnq(\ep)=\{G\in\ffnq\mid \ v_G(Q)=\ep\}$. Hence, for all $G\in\mathcal{F}(g)$, $G\not\in\ffnq(\ep)$, we have   $\ell\left(S_\ep(N^+_{\nu,Q}(G))\right)=0$.
	
	Indeed, if $G\not\in\ffnq$, then $Q\nmid_\nu G$ and Lemma \ref{length} shows that $\ell\left(N^+_{\nu,Q}(G)\right)=0$. If $G\in\ffnq$ but $v_G(Q)\ne\ep$ then $\ell\left(S_\ep(N^+_{\nu,Q}(G))\right)=0$ because $N^+_{\nu,Q}(G)$ is one-sided with a different slope.
	By Theorem \ref{product},  
	\begin{align*}
		\ell\left(S_\ep(N^+_{\nu,Q}(g))\right)=&\,\sum\nolimits_{G\in\mathcal{F}(g)}\ell\left(S_\ep(N^+_{\nu,Q}(G))\right)=\sum\nolimits_{G\in\ffnq(\ep)}\ell\left(S_\ep(N^+_{\nu,Q}(G))\right)\\=&\,\sum\nolimits_{G\in\ffnq(\ep)}
		(\deg(G)/\deg(Q).
	\end{align*}
	This proves (ii).
\end{proof}

\subsubsection{Newton polygons with respect to non-henselian valuations}

Let $(K,v)$ be an arbitrary  valued field and take a  monic $g\in \kx$. 
Let $\vb$ be a fixed extension of $v$ to $\kb$ and $(\kh,\vh)$ the henselization of $(K,v)$ determined by this choice.

Consider the unique extension $\muh$ of $\mu$ to $\khx$ whose restriction to $\kh$  is $\vh$ (Theorem \ref{Rig}). The strategy to prove Theorem \ref{conjNewton} is to deduce it from Theorem \ref{ConjHensel} after a suitable comparison of the sets $\ffmph$, $\ffmph(\la)$ with the analogous objects $\ffnq$, $\ffnq(\ep)$, with respect to $\nu=\muh$ and $Q=Q_{\mu,\phi}$.

To this end, we need a relevant consequence of Theorem \ref{NN}.

\begin{proposition}\cite{NN}\label{SpecialFactor}
	Let $\mu$ be a valuation on $\kx$ and let $\phi\in\kpm$ be a key polynomial of minimal degree. Let $Q=Q_{\mu,\phi}\in\irr(\kh)$ be the irreducible factor of $\phi$  determined by $\mu$. Then, $Q\in\kp(\muh)$ and \  $\inmh (\phi/Q)$ is a unit in $\ggmh$. 
\end{proposition}

Take a type $(\mu,\phi)$ on $\kx$ and let $Q\in\irr(\kh)$ be the irreducible factor of $\phi$ determined by $\mu$. From now on, we denote
$$
P=\phi/Q\in\khx,\qquad \al=\muh(P).
$$

The first thing to observe is that the types $(\mu,\phi)$ and $(\muh,Q)$ ``point out" to the same irreducible factors of $g$ in $\khx$.

\begin{lemma}\label{compF}
	$\ffmph=\mathcal{F}_{\muh,Q}(g)$.
\end{lemma}

\begin{proof}
	Let $\fff$ be the set of monic irreducible factors of $g$ in $\khx$. 
	
	By Theorem \ref{Rig}, for all $G\in\fff$, we have
	$$
	\mu<w_G\ \sii\ \muh<(w_G)^h=v_G.
	$$
	Thus, in order to prove the lemma, we need only to check that
	$$
	\ty(\mu,w_G)=[\phi]_\mu\ \sii\ \ty(\muh,v_G)=[Q]_{\muh}.
	$$
	By Corollary \ref{td}, this is equivalent to:
	$$
	\mu(\phi)<w_G(\phi)\ \sii\ \muh(Q)<v_G(Q).
	$$
	On the other hand,
	$$
	\muh(Q)=\muh(\phi)-\al=\mu(\phi)-\al,\quad
	v_G(Q)=v_G(\phi)-v_G(P)=w_G(\phi)-v_G(P).
	$$
	Thus, we must proof $\muh(P)=v_G(P)$. By Proposition \ref{SpecialFactor}, $\inmh P$ is a unit in $\ggmh$. Hence, $Q\nmid_{\muh}P$. Since $\ty(\muh,v_G)=[Q]_{\muh}$,  this implies $\muh(P)=v_G(P)$.
\end{proof}

\begin{corollary}\label{Thmi}
	Item (i) of Theorem \ref{conjNewton} follows from item (i) of Theorem \ref{ConjHensel}. 
\end{corollary}

\begin{corollary}\label{lalpha}
	For all $\la\in\gq$, $\ffmph(\la)=\mathcal{F}_{\muh,Q}(g)(\la-\al)
	$. 
\end{corollary}

\begin{proof}
	By Lemma \ref{compF}, we may rewrite these subsets as:
	$$
	\ars{1.3}
	\begin{array}{l}
	\ffmph(\la)=\{G\in\ffmph\mid w_G(\phi)=\la\},\\
	\mathcal{F}_{\muh,Q}(g)(\la-\al)=\{G\in\ffmph\mid v_G(Q)=\la-\al\}.
	\end{array}
	$$
	Along the proof of Lemma \ref{compF} we saw that $v_G(P)=\muh(P)=\al$. 
	Therefore, $w_G(\phi)=\la$ if and only if $v_G(Q)=\la-\al$.
\end{proof}\e

In order to finish the proof of Theorem \ref{conjNewton}, we must compare the Newton polygons $\npp(g)$ and $N^+_{\muh,Q}(g)$. Instead of comparing these polygons directly, we compare each one with an auxiliary polygon.
Consider the canonical $\phi$-expansion of $g$: 
$$
g(x)=\sum_{n\ge0} a_n\phi^n,\quad a_n\in\kx,\quad \deg(a_n)< \deg(\phi). 
$$

We may deduce from the $\phi$-expansion a trivial $Q$-expansion of $g$:
$$
g(x)=\sum_{n\ge0} b_nQ^n,\quad b_n=a_nP^n\in\khx. 
$$
This $Q$-expansion is far from being the canonical one, but it leads to Newton polygons which are easily comparable with $\np(g)$.\e

\nn{\bf Notation. }Let us denote by $\nnn(g)$ the convex hull of the cloud of points $$\{(n,\muh(b_n))\mid n\ge0\}.$$
Let $\nnn^+(g)$ be the polygon formed by the sides of $\nnn(g)$ of slope $-\ep$, with $\ep>\nu(Q)$.

\begin{lemma}\label{Comp1}
	The linear automorphism
	$$
	\Q\times\gq\lra \Q\times\gq,\qquad (x,y)\ \longmapsto\ (x,y+\al x)
	$$
	maps $\np(g)$ to $\nnn(g)$. Moreover, it maps $S_\la\left(\np(g)\right)$ to $S_{\la-\al}\left(\nnn(g)\right)$ and it preserves the lengths of these components. 
	In particular, it maps $\npp(g)$ to $\nnn^+(g)$.
\end{lemma}

\begin{proof}
	Since $\muh(b_n)=\mu(a_n)+n\alpha$, this linear automorphism maps:
	$$
	\{(n,\mu(a_n))\mid n\ge0\}\ \longmapsto\ \left\{(n,\muh(b_n))\mid n\ge0\right\}.
	$$
	Since linear mappings preserve convex subsets, $\np(g)$ is mapped to $\nnn(g)$.
	
	The points lying on $S_\la\left(\np(g)\right)$ correspond to monomials $a_n\phi^n$ satisfying:
	$$
	\mu(a_n)+n\la\le\mu(a_m)+m\la\quad\mbox{ for all }m\ge0.
	$$
	This is equivalent to 
	$$
	\muh(b_n)+n(\la-\al)\le\muh(b_m)+m(\la-\al)\quad\mbox{ for all }m\ge0.
	$$
	Thus, the linear automorphism maps $S_\la\left(\np(g)\right)$ to $S_{\la-\al}\left(\nnn(g)\right)$, and it preserves the lengths of these segments. 
\end{proof}\e

Therefore, the proof of Theorem \ref{conjNewton} follows immediately from Theorem \ref{ConjHensel}, once we prove the next result.

\begin{lemma}\label{Comp2}
	$\ \nnn^+(g)=N^+_{\muh,Q}(g)$.
\end{lemma}

\begin{proof}
	It suffices to show that, for all $\ep\in\gq$, $\ep>\muh(Q)$, we have
	$$
	S_\ep\left(\nnn^+(g)\right)=S_\ep\left(N^+_{\muh,Q}(g)\right).
	$$
	Consider the augmented valuation $\nep=[\muh; \,Q,\ep]$. Recall that $Q$ becomes a key polynomial of minimal degree for $\nep$.
	
	The monomials $b_nQ^n$ such that $\nep(b_nQ^n)=\nep(g)$ correspond to points lying on $S_\ep\left(\nnn^+(g)\right)$. Imagine that these monomials are:
	$$
	b_sQ^s+\cdots +b_tQ^t.
	$$
	Then, if we denote for simplicity $\bbr_n=\inn_{\nep} b_n$ and $\pi=\inn_{\nep}Q$, we have:
$$
	\inn_{\nep}(g)=\bbr_s\pi^s+\cdots +\bbr_t\pi^t\in\gg^0_{\nep}[\pi]=\gg_{\nep},\quad \bbr_n\in\gg^0_{\nep}.
$$
	Indeed, by Proposition \ref{SpecialFactor}, $\inmh b_n=\inmh(a_nP^n)$ is a unit in $\ggmh$. Hence, $Q\nmid_{\muh}b_n$, so that $\nep(b_n)=\muh(b_n)$. Hence, the homomorphism $\ggmh\to\gg_{\nep}$ maps $\inmh b_n$ to  $\bbr_n$, and the latter is a unit in $\gg_{\nep}$.
	
	Now, let $g=\sum_{n\ge0}c_nQ^n$ be the canonical $Q$-expansion of $g$. We can argue as above. The monomials $c_nQ^n$ such that $\nep(c_nQ^n)=\nep(g)$ correspond to points lying on $S_\ep\left(N^+_{\muh,Q}(g)\right)$. If these monomials are $	c_kQ^k+\cdots +c_\ell Q^\ell$,
	we deduce as above:
$$
	\inn_{\nep}(g)=\cb_k\pi^k+\cdots +\cb_\ell \pi^\ell\in\gg^0_{\nep}[\pi]=\gg_{\nep},\quad \cb_n\in\gg^0_{\nep}.
$$
	
	By Theorem \ref{g0gm}, we deduce that 
	$$
	s=k,\quad t=\ell,\quad \bbr_n=\cb_n  
	$$
	for all $s\le n\le t$ such that $\nep(b_nQ^n)=\nep(g)$, which must be the same indices for which   $\nep(c_nQ^n)=\nep(g)$.
	This proves the lemma.
\end{proof}

\section{The OM-algorithm}\label{secOM}
\subsection{A formal OM-algorithm}\label{subsecOM}\mbox{\null}\e

\nn{\bf OM-algorithm}

\nn{\bf Input: }$g\in\irr(K)$, \ $v$ a valuation on $K$

\nn{\bf Output: }A list of types singling out the irreducible factors of $g$ in $\khx$\e

\nn{\sl Stack}$\;\leftarrow\left[(v,x,\deg(g))\right]$; \
\nn{\sl Types}$\;\leftarrow[\ ]$

\nn{\bf while }$\#\mbox{\sl Stack}>0$ \ {\bf do}

\quad pick any $(\mu,\phi,\ell)\in$ {\sl Stack} and delete it from the {\sl Stack}

\quad{\bf for }$-\la$ slope of $\op{NP}(\mu,\phi,\ell)(g)$ {\bf do}

\quad\qquad $\mul\leftarrow [\mu;\,\phi,\la]$

\quad\qquad compute and factorize $R_{\mul,\phi}(g)=\psi_1^{n_1}\cdots \psi_s^{n_s}$ \;in \;$\kappa(\mul)[y]$  

\quad\qquad {\bf for } $1\le i\le s$  \ {\bf do}

\quad\qquad\qquad $\varphi\leftarrow\op{lift}_{\mul,\phi}(\psi_i)$

\quad\qquad\qquad {\bf if }$\deg(\varphi)>\deg(\phi)$ \,{\bf then }\,$\mu\leftarrow\mul$

\quad\qquad\qquad {\bf if }$n_i=1$ \,{\bf then }\,append $(\mul,\varphi)$ to {\sl Types}
{\bf else }append $(\mu,\varphi,n_i)$ to {\sl Stack} 

\nn{\bf return }{\sl Types}\bs

If the OM-algorithm terminates then it provides:

\begin{itemize}
\item An approximation in $\kx$ to each irreducible factor of $g$ in $\khx$.
\item All extensions of $v$ to the field $L=\kx/(g)$, plus a computation of their ramification indices and residual degrees.
\end{itemize}


Indeed, let  $(\mu,\phi)$ be an output type.
By storing all pairs $(\mu_n,\phi_{n+1})$ and slopes $\ga_{n+1}$ considered along the procedure, we get a MLV chain of ordinary augmentations:  
\begin{equation}\label{MLVwF}
v\ \stackrel{\phi_0,\ga_0}\lra\  \mu_0\ \stackrel{\phi_1,\ga_1}\lra\   \cdots
\ \lra\ \mu_{r-1} 
\ \stackrel{\phi_{r},\ga_{r}}\lra\ \mu_r=\mu  
\end{equation}
satisfying moreover $\deg(\phi_r)<\deg(\phi)$. 

\subsection*{Approximants to the irreducible factors}\label{subsubApp}
Let $G\in\irr(\khx)$ be the irreducible factor of $g$ singled out by the type $(\mu,\phi)$. By construction,  
$$
\ffmph=\left\{G\right\},\quad \ty(\mu,w_G)=[\phi]_\mu,\quad \deg(G)=\deg(Q_{\mu,\phi}).
$$

Since $\mu$ is inductive, Lemma \ref{IndRig} shows that $\phi\in\kpm\subset \kp(\muh)$ is irreducible over $\khx$. Thus,  $\phi=Q_{\mu,\phi}$ and $\deg(G)=\deg(\phi)$. 

We consider the whole class $[\phi]_\mu\subset \kpm$ as a set of  approximants to $G$ provided by the OM-algorithm. 

The measure of the quality of any $Q\in[\phi]_\mu$ as an approximation to $G$ is indicated by the value
$w_G(Q-G)=w_G(Q)$. This ``precision" is bounded above by
$$
\sup\left(w_G\left([\phi]_{\mu}\right)\right).
$$
If $v$ has rank one, then $K$ is dense in $\kh$ and this supremum is equal to infinity. In this case, there are approximants with arbitrarily large precision.  

For larger rank, this is not always possible. We give an example in Section \ref{subsubsecExample} where this supremum is bounded.

\subsection*{Extensions of $v$ to the field $L=\kx/(g)$}\label{subsubEF}

\begin{theorem}\label{OMesfs}
	Let $G\in\irr(\khx)$ be the irreducible factor of $g$ singled out by the type $(\mu,\phi)$. Then,
$$
e(\wb_G/v)=e_0\cdots e_r,\qquad f(\wb_G/v)=f_0\cdots f_r\deg(R_{\mu,\phi_r}(\phi)),
$$
where, $e_i, f_j$ are the numerical invariants of the MLV chain of $\mu$.
\end{theorem}

\begin{proof}
By Lemma \ref{IndRig}, $\phi\in\kp(\muh)$. Since $\ty(\mu,w_G)=[\phi]_\mu$, we have
$$
\muh(\phi)=\mu(\phi)<w_G(\phi)=v_G(\phi).
$$	
Hence, $\ty(\muh,v_G)=[\phi]_{\muh}$, by  Corollary \ref{td}. Now, since 
$\muh(G)<v_G(G)=\infty$,
we deduce that $\phi\mid_{\muh}G$. Since $\deg(G)=\deg(\phi)$, this implies that $G$ is a key polynomial for $\muh$ and $G\sim_{\muh} \phi$ \cite[Lem 2.5]{KP}. 

Therefore, it makes sense to consider the ordinary augmentation $[\muh;\, G,\infty]$. Since this valuation has support  $G\khx$, we have 
 $[\muh;\, G,\infty]=v_G$, by Proposition \ref{lfinh}.
By Lemma \ref{IndRig}, we obtain the following MLV chain of $v_G$:
$$
\vh\ \stackrel{\phi_0,\ga_0}\lra\  \muh_0\ \stackrel{\phi_1,\ga_1}\lra\   \cdots 
\ \stackrel{\phi_{r},\ga_{r}}\lra\ \muh_r=\muh  
\ \stackrel{G,\infty}\lra\  v_G,    
$$
whose initial numerical invariants $e_0,\dots, e_r;\,f_0,\dots,f_r$ coincide with those determined by the MLV (\ref{MLVwF}) of $\mu$.
By Proposition \ref{infinity}, 
$$
e(\vb_G/\vh)=e_0\cdots e_r,\qquad f(\vb_G/\vh)=f_0\cdots f_r\deg(R_{\muh,\phi_r}(G)).
$$

Since $G\sim_{\muh} \phi$, \cite[Cor. 5.5]{KP} shows that $R_{\muh,\phi_r}(G)=R_{\muh,\phi_r}(\phi)$. 
Also, let $u\in\ggm$ be a homogeneous unit such that  $\gr_\mu(u)=\mu(\phi_r^{e_r})$. Let $u^h\in\ggmh$ be the image of $u$ under the isomorphism $\ggm\hk\ggmh$. Then, it is easy to check that the following diagram commutes 
$$
\ars{1.3}
\begin{array}{rcl}
\kx&\subset&\khx\\
\raise.4ex\hbox{\mbox{\tiny $R_{\mu,\phi,u}$}}\downarrow\ \;&&\ \ \downarrow \raise.4ex\hbox{\mbox{\tiny $R_{\muh,\phi,u^h}$}}\\
\ka(\mu)[y]&\hra&\ka(\muh)[y]
\end{array}
$$
where the lower horizontal map is the isomorphism induced by $\ggm\hk\ggmh$. We deduce that $\deg\left(R_{\muh,\phi_r}(\phi)\right)=\deg\left(R_{\mu,\phi_r}(\phi)\right)$.

Therefore, the theorem will be proved if we check that
\begin{equation}\label{eeff}
e(\wb_G/v)=e(\vb_G/\vh),\qquad 
f(\wb_G/v)=f(\vb_G/\vh).
\end{equation}

Now, the extension $L^h=L\cdot \kh=\khx/(G)$
is a henselization of $(L,\wb_G)$. Hence, the commutative diagram of extensions of valuations: 
$$
\ars{1.2}
\begin{array}{rcl}
(\kh,\vh)&\lra&(L^h,\vb_G)\\
\uparrow\quad\;&&\quad\;\uparrow\\
(K,v)\ &\lra&(L,\wb_G)
\end{array}
$$
implies (\ref{eeff}), because $e(\vh/v)=1=f(\vh/v)$ and 
$e(\vb_G/\wb_G)=1=f(\vb_G/\wb_G)$.
\end{proof}\e

The piece of a MLV chain joining $\mu$ with $w_G$ is not easy to describe. It can contain several limit augmentations, as illustrated by the following example. 

\subsubsection{An example}\label{subsubsecExample}
Take a prime number $p\equiv1\md4$ and let $\ord_p$ be the $p$-adic valuation. Denote by $\bar{a}$ the reduction modulo $p$ of an integer $a\in\Z$.

Choose a $p$-adic root $i\in\Z_p$  of the polynomial $g=x^2+1$:
$$
i=i_0+i_1p^{\ell_1}+\cdots+ i_np^{\ell_n}+\cdots,
$$
with $0<i_n<p$ for all $n$. Denote the truncations of $i$ by
$$
a_n=i_0+i_1p^{\ell_1}+\cdots+ i_{n-1}p^{\ell_{n-1}}\in\Z.
$$

Consider the field $K=\Q(t)$ equipped with the $\ord_t$ valuation. Every $u\in K^*$ has an initial term $
\inn(u)=\left(u\,t^{-\ord_t(u)}\right)(0)\in \Q^*$, 
with respect to $\ord_t$.

Consider the following discrete rank-two valuation on $K$:
$$
v\colon K^*\lra \Z^2_{\lx},\qquad v(u)=\left(\ord_t(u),\ord_p(\inn(u))\right),
$$ 
with values in $\Z^2$ equipped with the lexicographical order.
The residue field is $k=\F_p$.

The OM-algorithm applied to the polynomial $g=x^2+1$ terminates after a single double-dissection. The Newton polygon $N_{v,x}(g)$ is one-sided of slope $(0,0)$ and for $\mu_0=[v;\,x,(0,0)]$, we have $R_{\mu_0,x}(g)=y^2+1=(y-\overline{i}_0)(y+\overline{i}_0)$.

Hence, the algorithm detects the irreducible factors  $x-i,\,x+i\in\khx$, singled out by the output types $(\mu_0,x-i_0)$, $(\mu_0,x+i_0)$, respectively  

We obtain MLV chains consisting of a single limit augmentation
$$
v\, \stackrel{x,(0,0)}\lra\, \mu_0 \, \lra\, w_{x-i}=[\cc,\,g,\infty],\qquad v\, \stackrel{x,(0,0)}\lra\, \mu_0 \, \lra\ w_{x+i}=[\cc',\,g,\infty], 
$$ 
where $\cc=\left([v;\,x-a_n,(0,\ell_n)]\right)_{n\ge0}$ and $\cc'=\left([v;\,x+a_n,(0,\ell_n)]\right)_{n\ge0}$. 

Note that the quality of the approximations is bounded:
$$
w_G\left([x-i_0]_{\mu_0}\right)=
\left\{\vb(a_n-i)\mid n\ge0\right\}\subset\{0\}\times \Z.
$$

\subsection{Termination of the OM-algorithm}

Let us first remark that all involved subroutines can be performed by real algorithms.\e

\nn{\bf Remark. }{\it
Suppose that $v$ has a finite rational rank ($\dim_\Q(\gq)<\infty$) and we have algorithms performing the following tasks.
\begin{itemize}
	\item Field operations in $K$ and  $k$ and computation of the valuation $v\colon K^* \twoheadrightarrow\g$.  
	\item Computation of the residue class $\oo_v^* \twoheadrightarrow k^*$ and a section $\op{lift}_v\colon k^*\to \oo_v^*$.  
	\item Polynomial factorization in $\kappa[y]$ for  arbitrary finite extensions $\kappa/k$. 
\end{itemize}	
Then, there are algorithms performing all subroutines of the  OM-algorithm. 
}\e

Indeed, we use only three subroutines:
$$
\op{NP}(\mu,\phi,\ell)(-),\qquad R_{\mul,\phi}(-),\qquad
\op{lift}_{\mul,\phi}(-).$$

The subroutine $\op{NP}(\mu,\phi,\ell)(-)$, described in Section \ref{subsecDissRPol}, requires only:

(i) \  A quotient-with-remainder routine in $\kx$ to compute truncated $\phi$-expansions.

(ii) \ A routine to compute $\mu$.

By Lemma \ref{ordinate}, the computation of $\mu(a)$ for any  $a\in\kx$ follows easily from the computation of the Newton polygon 
$$
\op{NP}\left(\mu_{r-1},\phi_{r-1},\lfloor\deg(a)/\deg(\phi_{r-1})\rfloor\right)(a).
$$
Thus, a recursive descending procedure along the MLV chain (\ref{MLVwF}), enables the computation of $\mu$, based in the end on the routine computing the valuation $v$. 

The routines $ R_{\mul,\phi}(-)$ and $\op{lift}_{\mul,\phi}(-)$ can be obtained by a similar descending recursive procedure, described in \cite[Sec. 5]{RPO}. \e

Therefore, the only obstacle for the termination of the OM-algorithm would be the existence of an infinite sequence of double-dissections (double {\bf for} loops).

Since $\fff$ is a finite set, it admits only a finite number of non-trivial dissections of any of its subsets.  Also, since $\deg(\phi)\le \deg(g)$ for all key polynomials $\phi$ constructed along the process, the condition $\deg(\phi)<\deg(\varphi)$ inside the second {\bf for} loop may occur only a finite number of times.  

Thus, the OM-algorithm does not terminate if and only if  there is an infinite sequence of  \emph{refinement steps}, defined as follows.
\e

\defn
A refinement step is a double {\bf for} loop which,  applied to a certain $(\mu,\phi,\ell)\in$\;{\sl Stack}, yields a unique triple $(\mul,\varphi,n)$, and moreover $\deg(\varphi)=\deg(\phi)$.\e

By Theorem \ref{charKP}, $\deg(\varphi)=\erel(\mul)\deg(\psi)\deg(\phi)$. Hence, a refinement step is characterized by the following two conditions:\e

$\bullet$ \ $\npp(g)$ is one-sided and its slope $-\la$ satisfies $\erel(\mul)=1$.

$\bullet$ \ $R_{\mul,\phi}(g)=(y-\zeta)^\ell$, for some  $\zeta\in\kappa(\mul)^*$.\e

In this case, we just replace $(\mu,\phi,\ell)$ with $(\mu,\varphi,\ell)$ in the {\sl Stack}. 

Montes proved that infinite sequences of refinement steps cannot occur in the discrete rank-one case \cite{montes}.

\begin{theorem}\label{classicOM}
If $v$ is discrete of rank-one, then the OM-algorithm terminates.
\end{theorem}	

Let us write $L=K(\t)$, where $\t\in L$ is the class of $x$ modulo the ideal $g\kx$.	

The proof of this theorem  is based on the finiteness of the local index 
$$
\op{ind}(g):= v\left((\oo_g\colon\oo_v[\t])\right)\in\g,
$$
where $\oo_g$ is the integral closure of $\oo_v$ in the finite extension $L/K$. 
Through an ordered isomorphism between $\g$ and $\Z$, this index  is identified with a non-negative integer. 
The theorem follows from the fact that in
all intermediate steps of the algorithm, including the refinement steps, there is a positive integer contributing to the total value of $\op{ind}(g)$ \cite[Thm. 4.8]{HN}.

\subsection{Infinite sequences of refinement steps}\label{sec:examples}
As we saw in Section \ref{subsecOM}, the OM-algorithm aims to construct a MLV chain 
$$v\ \stackrel{\phi_0,\ga_0}\lra\  \mu_0\ \stackrel{\phi_1,\ga_1}\lra\   \cdots
\ \stackrel{\phi_{r},\ga_{r}}\lra\ \mu_r=\mu  
$$
of a valuation $\mu$ which is ``sufficiently close" to the valuation $w_G$, for some  irreducible factor $G$ of $g$ in $\khx$. 

Denote $m_n=\deg(\mu_n)=\deg(\phi_n)$ for all $0\le n\le r$. For a field  $K\subset \mathbb{K}\subset \kh$, let
$$
V_{m_n}(\mathbb{K})=\left\{w_G(f)\mid f\in \mathbb{K}[x] \mbox{ monic},\ \deg(f)=m_n\right\}\subset \gq.
$$

For each $n>0$, the analysis of the augmentation step $\mu_{n-1}\to\mu_n$ leads to three different ``infinite refinement" situations:\e

(IR1) \ There exists $\max\left(V_{m_n}(K)\right)$,

(IR2) \ $\max\left(V_{m_n}(K)\right)$ does not exist, but there exists $\max\left(V_{m_n}(\kh)\right)$,

(IR3) \ $\max\left(V_{m_n}(\kh)\right)$ does not exist,\e

By   \cite[Thm. 4.7]
{MLV}, the augmentation $\mu_{n-1}\to\mu_n$ is ordinary in the case (IR1) and  a limit augmentation in cases (IR2) and (IR3). Vaqui\'e showed that limit augmentations in the henselian case occur only when $G$ has defect \cite{Vaq2}, \cite[Sec. 6]{MLV}. Thus, we say that  $\mu_{n-1}\to\mu_n$  is a \emph{defectless} limit augmentation  in the (IR2) case, and a \emph{defect} limit augmentation in the (IR3) case. Defect limit augmentations occur only when $\op{char}(k)=p>0$; also, $\deg(\mu_n)/\deg(\mu_{n-1})$ is necessarily a power of $p$.


\subsubsection{An example of (IR2)}\label{subsubsecIR2}

Let $(K,v)$ be the valued field considered in Section \ref{subsubsecExample}.
Let us apply the OM-algorithm to the polynomial 
$$
g=x^4+(t+2)x^2+1\in \kx.
$$  

The double-dissection applied to the triple $(v,x,4)$ yields a one-sided Newton polygon of slope $(0,0)$. For $\mu_0=[v,x,(0,0)]$, the residual polynomial factorizes
$$R_{\mu_0,x}(g)=1+2y^2+y^4=(y^2+1)^2=(y-\bar{i}_0)^2(y+\bar{i}_0)^2\in k[y].
$$
As lifts of the irreducible factors we may take $$\varphi=x-a_1=x-i_0, \qquad \varphi'=x+a_1=x+i_0.$$ We get {\sl Stack}\;$=\left[ (v,x-a_1,2),\,(v,x+a_1,2)\right]$. By Theorem \ref{conjNewton}, we detect a splitting of $g$ into a product of two (unknown) polynomials in $\khx$ of degree two.

The application of the double-dissection to the triple $(v,x-a_1,2)$ leads to an infinite sequence of refinements:
\begin{equation}\label{infSeq}
(v,x-a_1,2)\;\rightsquigarrow\; (v,x-a_2,2)\;\rightsquigarrow\quad \cdots\quad \rightsquigarrow \;(v,x-a_n,2)\;\rightsquigarrow\quad \cdots
\end{equation}
and a similar situation occurs for the triple $(v,x+a_1,2)$. 

Indeed, the truncated $(x-a_n)$-expansion of $g$ is $b_0+b_1(x-a_n)+b_2(x-a_n)^2$, with
$$
b_0=g(a_n)=c_n^2+ta_n^2,\quad b_1=g'(a_n)=4a_nc_n+2ta_n,\quad b_2=\dfrac12\,g''(a_n)=6c_n+t-4,
$$
where $c_n=a_n^2+1$. One checks easily that $\la_n:=v(c_n)=(0,\ell_n)$, so that  
$N^+_{v,x-a_n}(g)$ is one-sided of slope $-\la_n$ and contains the three points $(0,2\la_n)$, $(1,\la_n)$, $(2,(0,0))$.   

Denote $\rho_n=[v,x-a_n,\la_n]$. We may identify $\kappa(\rho_n)=k$ and take $u_n=\inn_{\rho_n}(p^{\ell_n})$ as a unit of grade $\la_n$. We get $R_{\rho_n,x-a_n}(g)=\left(y-\bar{i}_n\right)^2$ and a natural lift of $y-\bar{i}_n$ is $(x-a_n)-i_n p^{\ell_n}=x-a_{n+1}\in\kp(\rho_n)$.

The totally ordered family of valuations $\cc=\left(\rho_n\right)_n\ge0$ is an essential continuous family of augmentations of $\rho_0=\mu_0$. It can be easily shown that all polynomials of degree one are $\cc$-stable, but $\phi:=x^2+1$ is $\cc$-unstable:
$$
\phi=(x-a_n)^2+2a_n(x-a_n)+c_n\imp \rho_n(\phi)=\min\{2\la_n,\la_n,\la_n\}=\la_n,
$$
for all $n$. Thus, $\phi$ is a limit key polynomial for $\cc$ and the right triple to append to the {\sl Stack} would be $(\cc,\phi,2)$.
The double-precission loop can be applied to this triple, to continue the OM-algorithm.

The truncated $\phi$-expansion of degree two is the whole $\phi$-expansion:
$$
g=-t+t\phi+\phi^2.
$$  
The Newton polygon, displayed in Figure \ref{figOM}, is one-sided of slope $-\la=-(1,0)/2$.
 The limit augmentation $\mu=[\cc;\,\phi,\la]$ has $\erel(\mu)=2$. We may still identify $\kappa(\mu)=k$ and take $u=\inm t$ as a unit of grade $2\la$. We get $R_{\mu,\phi}(g)=y-1$. Thus, $g$ is a lift of $y-1$ and the algorithm appends the type $(\mu,g)$ to the output list {\sl Types}. 

We may proceed in a analogous way with the triple $(\cc',\phi,\la)$, where $\cc'=\left(\rho'_n\right)_{n\ge0}$ is the essential continuous family of the valuations $\rho'_n=[v,x+a_n,\la_n]$. 

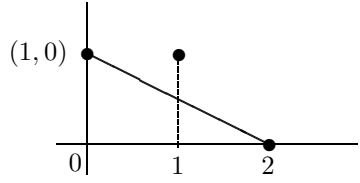
\begin{figure}
	\caption{Newton polygon $N_{\cc,x^2+1}(g)$ (or $N^+_{\mu_0,x-i}(g)$). }\label{figOM}
	\begin{center}
		\setlength{\unitlength}{4mm}
		\begin{picture}(10,5)
		\put(-.25,2.75){$\bullet$}\put(2.8,2.7){$\bullet$}\put(5.8,-0.3){$\bullet$}
		\put(-1,0){\line(1,0){10}}\put(0,-1){\line(0,1){5.7}}
		\put(0,3){\line(2,-1){6}}\put(0,3.03){\line(2,-1){6}}
		\put(2.8,-1){\begin{footnotesize}$1$\end{footnotesize}}
		\put(5.8,-1){\begin{footnotesize}$2$\end{footnotesize}}
		\put(-2.6,2.8){\begin{footnotesize}$(1,0)$\end{footnotesize}}
		\multiput(3,-.1)(0,.25){13}{\vrule height2pt}
		
		\put(-.6,-0.9){\begin{footnotesize}$0$\end{footnotesize}}
		\end{picture}
	\end{center}
\end{figure}\e

The OM-algorithm would output two types $(\mu,g)$, $(\mu',g)$ which single out two irreducible factors $G$, $G'$ of degree two of $g$ in $\khx$, with MLV chains: 
$$
v\ \stackrel{x,(0,0)}\lra\ \mu_0\ \stackrel{\phi,\la}\lra\ \mu\ \stackrel{g,\infty}\lra\ w_G,\qquad
v\ \stackrel{x,(0,0)}\lra\ \mu_0\ \stackrel{\phi,\la}\lra\ \mu'\ \stackrel{g,\infty}\lra\ w_{G'},
$$
where $\mu_0\to\mu$, $\mu_0\to\mu'$ are limit augmentations and  $\mu\to w_G$, $\mu'\to w_{G'}$ are ordinary augmentations. Since $\mu$ and $\mu'$ are not inductive, we do not get concrete approximations to the irreducible factors $G$, $G'$. However, we know their degrees, ramification indices and residual degrees.   Indeed, by Proposition \ref{infinity},
\begin{equation}\label{ef}
e(\bar{\mu}/v)=e(\bar{\mu}'/v)=2,\qquad 
f(\bar{\mu}/v)=f(\bar{\mu}'/v)=1.
\end{equation}

\subsubsection{An example of (IR1)}
Take $p$, $i$ as in the preceding example and consider the base field $K=\Q_p(t)$, equipped with the analogous discrete rank-two valuation 
$$
v\colon K^*\lra \Z^2_{\lx},\qquad v(u)=\left(\ord_t(u),\ord_p(\inn(u))\right).
$$ 

Let us apply the OM-algorithm to the same polynomial $g=x^4+(2+t)x^2+1$.

The double-dissection applied to the triple $(v,x,4)$ yields two triples $(v,x-a_1,2)$, $(v,x+a_1,2)$, each one leading to infinite sequences of refinements as in  (\ref{infSeq}).

However, imagine that our lifting routine chooses $$\varphi=\op{lift}_{\mu_0,x}(y-\bar{i}_0)=x-i,\qquad\varphi'=\op{lift}_{\mu_0,x}(y+\bar{i}_0)=x+i.$$ 
Then, the double-dissection applied to the triple $(v,x-i,2)$ is no more a refinement step. The truncated $(x-i)$-expansion of $g$ is $b_0+b_1(x-i)+b_2(x-i)^2$, with
$$
b_0=g(i)=-t,\quad b_1=g'(i)=2ti,\quad b_2=g''(i)/2=t-4.
$$
Thus, $N^+_{\mu_0,x-i}(g)$ is the  polygon displayed in Figure \ref{figOM}. It is one-sided of slope $-\la=-(1,0)/2$. 
The ordinary augmentation $\mu=[\mu_0;\,x-i,\la]$ has $\erel(\mu)=2$. We may still identify $\kappa(\mu)=k$ and take $u=\inm t$ as a unit of grade $2\la$. We get
$$
\inm g=\inm(-t-4(x-i)^2)=
-4u\left(\frac14+\frac{(x-i)^2}{u}\right),
$$
so that $R_{\mu,x-i}(g)=y+(1/4)$, admitting $\phi=(x-i)^2+(t/4)$ as a lift. The algorithm appends the type $(\mu,\phi)$ to the output list {\sl Types}. 

We may proceed in a analogous way with the triple $(v,x+i,2)$ to obtain the augmentation $\mu'=[\mu_0;\,x+i,\la]$ and a key polynomial $\phi'=(x+i)^2+(t/4)$.
The output of the OM-algorithm is a list of two types $(\mu,\phi)$, $(\mu',\phi')$ which single out two irreducible factors $G$, $G'$ of degree two of $g$ in $\khx$. 

The MLV chains of $\mu$ and $\mu'$ contain only ordinary augmentations:
$$
v\ \stackrel{x,(0,0)}\lra\ \mu_0\ \stackrel{x-i,\la}\lra\ \mu,\qquad
v\ \stackrel{x,(0,0)}\lra\ \mu_0\ \stackrel{x+i,\la}\lra\ \mu'
$$
Hence, the key polynomials $\phi=(x-i)^2+(t/4)$, $\phi=(x+i)^2+(t/4)$ are approximations to the true factors $G$, $G'$, respectively.   
We have
$$
w_G(x-i)=\la=(1/2,0)=\max\left(V_1(K)\right).
$$ 

The ramification indices and residual degrees of $G$ and $G'$ are given by (\ref{ef}).

\subsubsection{An example of (IR3)}

Let $\F$ be an algebraic closure of the prime field $\F_p$, for some prime number $p$. For an indeterminate $t$, consider the fields of Newton-Puiseux series and Hahn series in $t$, respectively:
$$
K=\bigcup_{N\in\N} \F((t^{1/N}))\subset \F((t^\Q)).
$$
The Hahn field $\F((t^\Q))$ consists of all power series with rational exponents and well-ordered support. For instance, as remarked by Abhyankar,
$$
s=\sum\nolimits_{n\ge1}t^{-1/p^n}\in \F((t^\Q)),
$$
is a root of the Artin-Schreier irreducible polynomial $g=x^p-x-t^{-1}\in \kx$. The truncations of $s$ belong to $K$:
$$
s_n=t^{-1/p}+\cdots+t^{-1/p^{n}}\in K,\quad n\ge 1.
$$ 

Consider the valuation $v=\ord_p$ on $K$, with value group  $\g=\Q$ and residue field $k=\F$. The valued field $(K,v)$ is henselian.

Let us apply the OM-algorithm to test the irreducibility of $g$.
The double-dissection applied to the triple $(v,x,p)$ yields a one-sided Newton polygon of slope $1/p$. For $\mu_0=[v,x,-1/p]$, and $u_0=\inn_{\mu_0}t^{-1/p}$ as a chosen unit of grade $-1/p$, the residual polynomial is $R_{\mu_0,x}(g)=(y-1)^p\in k[y]$.
Take $\varphi=x-s_1$ as a lift of $y-1$. 
The iterative application of the double-dissection  leads to an infinite sequence of refinements:
$$
(v,x-s_1,p)\;\rightsquigarrow\; (v,x-s_2,p)\;\rightsquigarrow\quad \cdots\quad \rightsquigarrow \;(v,x-s_n,p)\;\rightsquigarrow\quad \cdots
$$

Consider the essential continuous family $\cc=(\rho_n)_{n\ge0}$, where $\rho_n=[v;\ x-s_n,1/p^{n+1}]$.
All polygons $N_{v,x-s_n}(g)$ are one-sided of slope $1/p^{n+1}$ and  $R_{\rho_n,x-s_n}(g)=(y-1)^p$, if we choose $u_n=\inn_{\mu_n}t^{-1/p^{n+1}}$ as a unit of grade $-1/p^{n+1}$.

The polynomial $g$ is a limit key polynomial, and  $w_g=[\cc;\,g,\infty]$. 
The unique extension of $v$ to $L=\kx/(g)$ is the valuation $\bar{w}$ naturally induced by $w_g$. It has 
$$
e(\bar{w}/v)=f(\bar{w}/v)=1, \qquad d(\bar{w}/v)=p,
$$
where $d(\bar{w}/v)$ is the defect of the extension.\e

\nn{\bf Conclusion. }{\it In order to overcome the existence of an infinite sequence of refinements of type (IR1), the lifting routine $\op{lift}_{\mu_{n-1},\phi_{n-1}}(-)$ should be clever enough to compute $\max(V_{m_n}(K))$ in a finite number of steps, for all $n$.

In order to overcome the existence of an infinite sequence of refinements of types (IR2) or (IR3), the OM-algorithm should be modified to enable it to detect limit augmentations and compute limit key polynomials}.

 \section{A polynomial factorization algorithm in the rank-one case}

In this section, we analyze two OM-based algorithms of Poteaux-Weimann: a polynomial irreducibility test and a polynomial factorization algorithm \cite{PW2}. Both algorithms assume that the base valuation $v$ is discrete of rank one and impose certain conditions on the input polynomial $g\in\kx$. 


We generalize the irreducibility test to an arbitrary valued field $(K,v)$. The condition imposed on $g$ ensures that no defect limit augmentation appears. Since defectless limit augmentations only appear for input polynomials which are not irreducible in $\khx$, they are no obstacle for an irreducibility test. The crucial feature of the agorithm is that infinite sequences of refinements of type (IR1)  are avoided by considering \emph{approximate roots} as ``optimal" key polynomials.
 

As a consequence, we obtain a polynomial factorization over $\khx$ for a (not necessarily discrete) base valued field $(K,v)$ of rank one. In this case, $K$ is dense in $\kh$ so that defectless limit augmentations cannot appear. Also, defect limit augmentations are avoided by a strong condition imposed on the input polynomial. 

\subsection{Irreducibility test}\label{subsecTest}

Let $g\in\kx$ be a monic polynomial. Let $n$ be a divisor of $\deg(g)$ such that $\chr(K)\nmid n$.

The \emph{approximate root} $Q=\root{n}\of{g}$ is a monic polynomial in $\kx$, of degree $\deg(g)/n$, such that the canonical $Q$-expansion of $g$:
$$
g=Q^n+a_{n-1}Q^{n-1}+\cdots +a_1Q+a_0, \quad \deg(a_i)<\deg(Q),
$$
satisfies $a_{n-1}=0$.

Approximate roots were introduced by Abhyankhar and Moh in \cite{AM} as a tool to prove the embedding line theorem (see \cite{Pop} for a survey). In \cite{Ab}, Abhyankhar used approximate roots for an irreducibility test in $\mathbb{C}[[x]][y]$, then generalized in \cite{PW1,PW2} over a complete discrete valuation ring.

It is obvious that the approximate root is unique, if it exists. The existence follows from the following result, which gives moreover a concrete algorithm to compute it.

\begin{lemma}\cite[proof of Proposition 6.3]{Pop}\label{existence}
Let $g\in\kx$ be a monic polynomial. Let $n$ be a divisor of $\deg(g)$ such that $\chr(K)\nmid n$.
Take any monic polynomial $\phi\in\kx$ of degree $\deg(g)/n$, and consider the $\phi$-expansion
$$
g=\phi^n+a_{n-1}\phi^{n-1}+\cdots +a_1\phi+a_0, \quad \deg(a_i)<\deg(\phi).
$$
Take $\phi^*=\phi+(a_{n-1}/n)$ and let $a^*_{n-1}$ be the $(n-1)$-th coefficient of the $\phi^*$-expansion of $g$. Then, if $a_{n-1}\ne0$, we have $\deg(a^*_{n-1})<\deg(a_{n-1})$.
\end{lemma}

The next result establishes a link between approximate roots and key polynomials. 

\begin{proposition}\label{AppKP}
Let $\mu$ be a valuation on $\kx$ and $\varphi$ a key polynomial for $\mu$.
Let $g\in\kx$ be a monic polynomial such that $\chr(\km)\nmid \deg(g)$ and moreover: 
\begin{enumerate}
	\item[(i)\,] $N_{\mu,\varphi}(g)$ is one-sided of slope $-\la$.  
	\item[(ii)] For $\mul=[\mu;\,\varphi,\la]$, we have $R_{\mul,\varphi}(g)=\psi^n$, for some $\psi\in\irr(\ka(\mul))$.
\end{enumerate}

Then, the following holds:
\begin{enumerate}
	\item[(a)] The approximate root $Q=\root{n}\of g$ is a key polynomial for $\mul$ and $R_{\mul,\varphi}(Q)=\psi$.  
	\item[(b)] If $N_{\mul,Q}(g)$ is one-sided of slope $-\la_*$ and for $\nu=[\mul;\,Q,\la_*]$ we have $R_{\nu,Q}(g)=\psi_*^{n_*}$  for some $\psi_*\in\irr(\ka(\nu))$, then $n_*<n$. 
\end{enumerate}
\end{proposition} 

\begin{proof}
By (i) and (ii), $\deg(g)=efn\deg(\varphi)$, where 
$e=\erel(\mul)$ and $f=\deg(\psi)$.

Take any monic $\phi\in\kx$ of degree $ef\deg(\varphi)=\deg(Q)$ such that $R_{\mul,\varphi}(\phi)=\psi$. By Theorem \ref{charKP}, $\phi$ is a key polynomial for $\mul$. 

By Lemma \ref{existence}, we may obtain $Q$ from $\phi$ by a finite number of transformations of the form $\phi\mapsto \phi^*=\phi+(a_{n-1}/n)$, where $a_{n-1}$ is the $(n-1)$-th coefficient of the $\phi$-expansion of $g$. Hence, in order to prove (a), it suffices to show that $\phi^*\in\kp(\mul)$ and 
$R_{\mul,\varphi}(\phi^*)=\psi$.

By Lemma \ref{length}, the Newton polygon $N^+_{\mul,\phi}(g)$ has length $n$. Since $$n=\deg(g)/\deg(\phi)=\ell\left(N_{\mul,\phi}(g)\right),$$ we deduce that $N_{\mul,\phi}(g)=N^+_{\mul,\phi}(g)$. Thus, all slopes $-\ep$ of $N_{\mul,\phi}(g)$ satisfy $\ep>\mul(\phi)$. 
Now,  the point of abscissa $n-1$ lying on  $N_{\mul,\phi}(g)$ is $(n-1,\ep)$, where $-\ep$ is the largest slope of this polygon. Since the point $(n-1,\mul(a_{n-1}))$ lies on or above the polygon, we have $\mul(a_{n-1})\ge\ep$.  Since 
$\chr(\km)\nmid \deg(g)$, we have $\mul(n)=0$, so that
$$
\mul(a_{n-1}/n)=\mul(a_{n-1})\ge\ep>\mul(\phi).
$$
Thus, $\phi\sim_{\mul}\phi^*$. Since $\deg(\phi)=\deg(\phi^*)$, we deduce that $\phi^*$ is a key polynomial for $\mul$ and $R_{\mul,\varphi}(\phi^*)=R_{\mul,\varphi}(\phi)=\psi$ \cite[Lem. 2.5, Cor. 5.5]{KP}. This proves (a).

Let $g=\sum_{i\ge0}b_iQ^i$ be the $Q$-expansion of $g$. By the definition of the approximate root, $b_{n-1}=0$.
Under the hypotheses of (b), we have $$n\deg(Q)=\deg(g)=e_*f_*n_*\deg(Q),$$ 
where $e_*=\erel(\nu)$ is the least positive integer such that $e_*\la_*\in\g^0_{\nu}=\g_{\mul}$ and $f_*=\deg(\psi_*)$. Hence, $n=e_*f_* n_*$ and the equality $n=n_*$ holds only when $e_*=f_*=1$. This is incompatible with the assumptions in (b). Indeed, suppose that $n=n_*$. Since  $f_*=1$ we have $R_{\nu,Q}(g)=(y+\zeta)^n$ for some $\zeta\in\ka(\nu)^*$. Since $\chr(\mul)\nmid n$, the $(n-1)$-th coefficient of this polynomial is $n\zeta^{n-1}\ne0$. On the other hand, since $e_*=1$, by the definition of the residual coefficients in (\ref{resCoeff}), the point $(n-1,b_{n-1})$ lies on the Newton polygon. This is a contradiction because $b_{n-1}=0$. This proves (b).
\end{proof}\e

Therefore, the following irreducibility test works for arbitrary valued fields $(K,v)$.\e

\nn{\bf Irreducibility test}

\nn{\bf Input: }$(K,v)$ valued field, $g\in\kx$ monic, square-free such that $\chr(k)\nmid\deg(g)$

\nn{\bf Output: }A boolean (is $g$ irreducible over $\khx$?)\e

\nn$\mu\leftarrow v$; $\;\phi\leftarrow x$; $\;n\leftarrow \deg(g)$

\nn{\bf while }$n>1$ \ {\bf do}

\quad {\bf if } $N_{\mu,\phi}(g)$ is one-sided (of slope $-\la$) \ {\bf then } $\mu\leftarrow [\mu;\,\phi,\la]$ {\bf \ else }
{\bf return } False

\quad {\bf if } $R_{\mu,\phi}(g)=\psi^m$ for some $\psi\in\irr(\kam)$ {\bf \ then } $\phi\leftarrow \root{m}\of{g}$ {\bf \ else }  {\bf return } False

\quad$n\leftarrow m$

\nn{\bf return } True\e

In the very first step, it could happen that $N_{v,x}(g)$ is one-sided of slope $-\la\in\g$ and $R_{\mu,x}(g)$ is the $n$-th power of a polynomial of degree one. With the notation of Proposition \ref{AppKP}, we would have $e=f=1$ and $n=\deg(g)$. However, Proposition \ref{AppKP} shows that, as long as we do not detect a factorization of $g$ over $\khx$, in all further steps we will have $e_*f_*>1$, so that no refinement steps occur. Therefore, the algorithm terminates in $\log(\deg(g))$ steps.

As we did for the OM-algorithm, if $g$ is irreducible over $\khx$, then by storing all types $(\mu,\phi)$ obtained along the process,  we obtain a MLV chain of its associated valuation $w_g$, plus a computation of the residual degree and ramification index of the unique extension of $v$ to the field $\kx/(g)$.

 \subsection{Polynomial factorization}\label{subsecFactor}

Consider a monic, square-free $g\in \kx$ such that $\chr(k)\nmid \deg(g)$.

A \emph{splitting pair} of $g$ is any pair $(\mu,\phi)$ considered in the last call of the {\bf while} loop of the Irreducibility test of Section \ref{subsecTest}.

Note that, either $(\mu,\phi)=(v,x)$, or $(\mu,\phi)$ is a type.   
A splitting pair has the following general properties.

\begin{lemma}\label{lem:splitpair}
Let $g\in K[x]$ monic square-free with splitting pair  $(\mu,\phi)$. 
\begin{enumerate}
\item[(i)] If $(\mu,\phi)\ne(v,x)$, then $\ffmph=\ff(g)$.
\item[(ii)] $\np(g)=\npp(g)$ 
\item[(iii)] $\deg(g)=n\deg(\phi)$, where $n=\ell(N_{\mu,\phi}(g))$.
\item[(iv)] $\phi$ is irreducible over $\khx$.  
\end{enumerate}
\end{lemma}

\begin{proof}
Items (i), (ii), (iii) follow immediately from the design of the Irreducibility test. Item (iv) follows from Lemma \ref{IndRig}, because $\mu$ is an inductive valuation. 
\end{proof}

\e

\subsubsection{Right end-slope factorization}
Let $-\lambda$ the \emph{right end-slope} of $N=\np(g)$. That is, the slopes $-\ep$ of $N$ satisfy $-\ep\le-\la$, or equivalently, $\ep\ge\la$. Let $n_\la$ be the abscissa of the left end-point of $S_\la(N)$ (cf. Figure \ref{figComponent1}).

From now on, we assume $g$ reducible in $K^h[x]$ and denote $\mul=[\mu,\phi,\la]$.

Although a splitting pair of $g$ is not unique, this valuation $\mul$ is intrinsically associated to $g$. Indeed, Figure \ref{figDD} shows that $\mul$ is the greatest common lower node in $\ttt$ of the finite set of leaves $\{w_G\mid G\in\ff(g)\}$. The existence of greatest common lower nodes in $\ttt$ is guaranteed by \cite[Prop. 5.2]{VT}.

By the definition of a splitting pair, at least one of the following situations occurs:
\begin{itemize}
\item $N_{\mu,\phi}(g)$ is not one-sided; that is, $n_\la>0$.
\item  $R_{\mul,\phi}(g)=\psi_1^{n_1}\cdots\psi_s^{n_s}$, \,with $\psi_i\in\irr(\kappa(\mul))$ and $s\ge 2$.
\end{itemize}
\e
Let $\varphi_i=\op{lift}_{\mul,\phi}(\psi_i)\in \kp(\mul)$ be monic lifts of $\psi_i$ for $i=1,\ldots,s$. We denote $n_0=n_\la$ and $\varphi_0=\phi$ for convenience. We have $n_0> 0$ or $s\ge 2$. 

\begin{lemma}\label{lem:initialisation}
We have $g\sim_{\mul} \prod_{i=0}^s \varphi_i^{n_i}$ and $\deg(g)=\sum_{i=0}^s \deg(\varphi_i^{n_i})$.
\end{lemma}

\begin{proof}
By \eqref{piFactors} in Section \ref{subsecResPol}, we know that $
\inn_{\mul}g\,\sim_{\unit}\,\inn_{\mul}(\prod_{i=0}^s \varphi_i^{n_i})$. By Lemma \ref{lem:splitpair}, and our choice of $\lambda$, we get $\lc_{\mul}(g)=1$. Since $\phi\in \kp(\mul)$ has minimal degree and $\varphi_i\in \kp(\mul)$ for all $i\ge 0$, Theorem \ref{charKP} implies that $\lc_{\mul}(\varphi_i)=1$.

By Theorem \ref{g0gm}, $\inn_{\mul}g=\inn_{\mul}(\prod_{i=0}^s \varphi_i^{n_i})$, because both elements have leading coefficient $1$.  

On the other hand, $\dgm(g)=n$, $\dgm(\varphi_0)=1$ and, again by Theorem \ref{charKP},  $\dgm(\varphi_i)=\deg(\varphi_i)/\deg(\phi)$ for all $1\le i\le s$.

By Theorem \ref{g0gm}, $n=\sum_{i=0}^s\dgm(\varphi_i)^{n_i}$. By multiplying both sides of this equality by $\deg(\phi)$, we deduce that $\deg(g)=\sum_{i=0}^s \deg(\varphi_i^{n_i})$.
\end{proof}

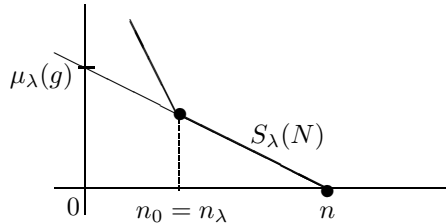
\begin{figure}
	\caption{Right end-side of $N=\np(g)$, defined by a splitting pair. We have $\mul(g)=n\la=\mul(\phi^n)$, for $n=\deg(g)/\deg(\phi)$.}\label{figComponent1}
	\begin{center}
		\setlength{\unitlength}{4mm}
		\begin{picture}(30,7)
		\put(10.9,2.7){$\bullet$}
		\put(15.8,0.15){$\bullet$}
		\put(7.8,4.5){\line(1,0){0.5}}
		\put(7,0.5){\line(1,0){13}}
		\put(8,-0.4){\line(0,1){7}}
		\put(7,5){\line(2,-1){9}}
		\put(11,3.1){\line(-1,2){1.5}}
		\put(11,3.04){\line(-1,2){1.5}}
		\put(11,3){\line(2,-1){5}}
		\put(11,3.04){\line(2,-1){5}}
		\multiput(11.1,.4)(0,.25){10}{\vrule height2pt}
		\put(15.8,-.5){\begin{footnotesize}$n$\end{footnotesize}}
		\put(9.7,-.5){\begin{footnotesize}$n_0=n_{\la}$\end{footnotesize}}
		\put(7.4,-.4){\begin{footnotesize}$0$\end{footnotesize}}
		\put(5.5,4){\begin{footnotesize}$\mul(g)$\end{footnotesize}}
		\put(13.5,2){\begin{footnotesize}$S_\la(N)$\end{footnotesize}}
		\end{picture}
	\end{center}
\end{figure}\e

Let $\mul^h$  be the unique common extension of $\mul$ and $\vh$ to $\khx$ (Theorem \ref{Rig}).

\begin{proposition}\label{prop:factorization}
There exist unique monic polynomials $G_0,\ldots,G_s\in K^h[x]$ such that 
$g=G_0\cdots G_s$, $G_i\sim_{\mul^h} \varphi_i^{n_i}$ and $\deg(G_i)=\deg(\varphi_i^{n_i})$ for all $i$. If $n_i=1$, then $G_i$ is irreducible. 
\end{proposition}

\begin{proof}
Since the augmentation $\mu\to \mul$ is ordinary, the valuation $\mul$ is inductive. Thus, $\varphi_0,\ldots,\varphi_s\in\kp(\mul^h)$, by Lemma \ref{IndRig}. 

For all $0\le i\le s$, let $G_i$ be the product of all irreducible factors $G$ of $g$ in $\khx$ satisfying $\varphi_i\mid_{\mul^h}G$. 
Then, the result follows from Lemma \ref{lem:initialisation} and Theorem \ref{fundamental}. 
\end{proof}

\subsubsection{Hensel lifting}
 The next result is a generalization of the multifactor Hensel lifting \cite[Algorithm 15.17]{Gathen} to an arbitrary valuation. We keep the notation of the previous paragraph. 

\begin{proposition}\label{propHensel}
Let $\gamma=\mul(g-\varphi_0^{n_0}\cdots \varphi_s^{n_s})-\mul(g)$. For all $n\in \N$ we can compute  monic polynomials $G_0^{(n)},\ldots,G_s^{(n)}\in K[x]$ such that $\mul^h(G_i-G_i^{(n)})>\mul^h(G_i)+2^n\gamma$.
\end{proposition}

\begin{proof}
Note that $\gamma > 0$ by Lemma \ref{lem:initialisation}.
Such a valuated Hensel lifting is detailed in \cite[Section 4.3]{PW2} in the discrete rank-one case. Since $\varphi_i^{n_i}$ is strongly monic in $\phi$ with respect to $\mul$ (\cite[Definition 5]{PW2}), then \cite[Lemma 7]{PW2} remains true in our context and \cite[Algorithm \texttt{HenselStep}]{PW2} extends straightforwardly to the valuation $\mul$.
\end{proof}

\e

This Hensel-like algorithm has quadratic convergence in the sense that the precision is doubled at each Hensel step. 
The following corollary is immediate. 

\begin{corollary}\label{cor:converge}
If the sequence $(2^n \gamma)_{n\in \N}$ is unbounded, then the sequence $(G_i^{(n)})_{n\in \N}$ converges to $G_i$ for all $0\le i\le s$. In particular, this holds whenever $v$ has rank one.
\end{corollary}

If $v$ has rank one, any choice of lifts $\varphi_i$ will allow to approximate the $G_i$'s with an arbitrary precision. If $v$ has rank $>1$, this is not always possible, as illustrated by the examples in Section \ref{sec:examples}. 

\subsubsection{Gauss valuation and Okutsu bound}\label{ssec:gaussprecision} In order to factorize recursively each approximant $G_i^{(n)}$ of Proposition \ref{propHensel}, we will rather measure the approximation with the Gauss valuation $v_0^h:\khx\to \Gamma$, $v_0^h(\sum_i c_i x^i):=\min v^h(c_i)$, which offers the advantage to be independent of the current splitting pair $(\mu,\phi)$. The valuation $v_0^h$ is asymptotically equivalent to the valuation $\mul^h$ in the following sense : 

\begin{lemma}\label{lem:compareValuations}
Suppose $g\in \oo[x]$ monic, and let $\mu,\phi,\lambda$ as above. Then, for all $F\in \khx$, we have $$v_0^h(F)\le\muh(F)\le \mul^h(F)\le v_0^h(F)+\lambda \frac{\deg(F)}{\deg(\phi)}.$$
\end{lemma}

\begin{proof}
The assumption $g\in \oo[x]$ ensures that the right end-slope $-\lambda_0$ of $N_{v,x}(g)$ satisfies $\lambda_0\ge 0$. Hence, the first extended valuation $\mu_0^h:=[v^h,x,\lambda_0]$ of the MLV chain of $\mul^h$ computed by the irreducibility test of $g$ satisfies $\mu_0^h\ge v_0^h$. Since $\mul^h\ge \mu_0^h$, this proves the left inequality of the lemma. Let us prove the right inequality. Since $\mul^h$ and $v_0^h$ coincide on $\kh$, we may suppose $F\in \oo^h[x]$ up to multiplying $F$ by a suitable constant $c\in \kh$. Also, it's enough to consider the case $F$ irreducible in $\oo^h[x]$. In such a case, $v_0^h(F)=0$ and the claim follows from \cite[Theorem 3.9]{KP}, having in mind that $\la=\mul(\phi)$. 
\end{proof}
\e

In what follows, splitting pairs of polynomials in $\kh[x]$ are defined as for $K[x]$.

\begin{definition}\label{def:okutsubound}
Let $F\in K^h[x]$ monic square-free, with splitting pair $(\mu^h,\phi^h)$ and right end-slope $-\lambda$. We define the \emph{Okutsu bound} of $F$ as 
$$
\delta_0(F):=\begin{cases} \mu^h(F), \quad {\rm if}\,\, F \,\,{\rm is\,\,irreducible} \\
\mul^h(F),\quad {\rm if}\,\, F\,\, {\rm is \,\,reducible} 
\end{cases}$$
\end{definition}

The notation and terminology for $\delta_0(F)$ is justified by the fact that Definition \ref{def:okutsubound} coincides with \cite[Definition 5.9]{defless} when $F$ is irreducible. 
In particular, $\delta_0(F)$ does not depend on the choice of the splitting pair.

\begin{proposition}\label{prop:precision}
Let $f,g\in \oo[x]$ be monic, square-free such that $v_0(f-g)> \delta_0(g)$ and $char(k)\nmid \deg(f)\deg(g)$. Then, $g$ is irreducible in $\khx$ if and only if $f$ is. 

If $g$ is reducible, then $g$ and $f$ have the same right end-slope and the same right end residual polynomial.
\end{proposition}

\begin{proof}
Since $g\in \oo[x]$, Lemma \ref{lem:compareValuations} shows that
\begin{equation}\label{lem68}
\mul(f-g)\ge\mu(f-g)\ge v_0(f-g)> \delta_0(g)\ge\mu(g)\ge v_0(g)\ge0.
\end{equation}
Hence,  $\deg(g)=\deg(f)$,  $g\sim_{\mu} f$ and $\np(g)=\np(f)$. If $g$ is irreducible, then $\np(f)=\np(g)$ has length $1$, so that $f$ is irreducible too. If $g$ is reducible, then $\delta_0(g)=\mul(g)$ and (\ref{lem68}) implies  $g\sim_{\mul} f$. Hence, $R_{\mul,\phi}(g)=R_{\mul,\phi}(f)$ \cite[Cor. 5.5]{KP} and $f$ is reducible too.  
\end{proof}

\e

\begin{corollary}\label{cor:precision}
Let $g\in \oo[x]$ monic square-free such that $char(k)\nmid \deg(g)$. Running algorithm \emph{Irreducibility}($g$) with Gauss precision $\sigma > \delta_0(g)$ returns a correct answer and allows to compute a splitting pair $(\mu,\phi)$ of $g$. If $g$ is reducible, the precision $\sigma$ is also sufficient to compute the right end-slope $\lambda$ and $R_{\mul,\phi}(g)$. 
\end{corollary}

\subsubsection{A factorization algorithm} 
For $g\in K^h[x]$ monic square-free, we define  $$\delta_{\max}(g):=\max\{\delta_0(G), \,\, G\in \fff\}.$$
Previous results lead to the following algorithm.
\e

\nn{\bf Factorization algorithm}

\nn{\bf Input: }$g\in \oo[x]$ monic square-free with $\chr(k)=0$ or $\chr(k)>\deg(g)$ and $\sigma\in \Gamma$ such that $\sigma>\delta_{\max}(g)$.

\nn{\bf Output: } The irreducible factors of $g$ in $\khx$ computed with Gauss precision $\ge \sigma$.

\e

\nn Run {\sl Irreducibility}($g$) with precision greater than $\delta_0(g)$.

\nn{\bf if} $g$ is irreducible {\bf then} {\bf return }$[g]$ {\bf else}:

\quad \nn $(\mu,\phi) \gets$ splitting pair of $g$

\quad  \nn $-\la\gets $ right end-slope of $\np(g)$

\quad \nn $\mul\leftarrow [\mu;\,\phi,\la]$

\quad \nn Compute and factorize $R_{\mul,\phi}(g)=\psi_1^{n_1}\cdots \psi_s^{n_s}$ \;in \;$\kappa(\mul)[y]$  

\quad \nn Compute some $\varphi_i\leftarrow\op{lift}_{\mul,\phi}(\psi_i)$ for $1\le i\le s$ and let $(\varphi_0,n_0)\gets (\phi,n_\lambda)$

\quad \nn $\gamma\gets \mul(g-\varphi_0^{n_0}\cdots \varphi_s^{n_s})-\mul(g)$. 

\quad \nn Compute $n\in \N$ such that $2^n\gamma\ge \sigma+\lambda \deg(g)/\deg(\phi)$

\quad \nn Compute $G_0^{(n)},\ldots,G_s^{(n)}\in \oo[x]$
as in Proposition \ref{propHensel}

\quad \nn {\sl Res} $\gets [\,]$

\quad \nn {\bf for } $i=0,\ldots,s$ {\bf do}:

\quad \quad {\bf if} $n_i=1$ \,{\bf then }\,append $G_i^{(n)}$ to {\sl Res} {\bf else } append {\sl Factorization}($G_i^{(n)}$) to {\sl Res}

\quad \nn{\bf return }{\sl Res}

\e

\begin{theorem}
If $v$ has rank one, then the algorithm {\sl Factorization} terminates and returns a correct answer. Moreover, the approximant factors converge to the irreducible factors of $g$ when we let $n\to+\infty$.
\end{theorem}

\begin{proof}
Since $\Gamma$ has rank one and $\gamma>0$, there exists $n\in\N$ such that $$2^n\gamma\ge \sigma+\lambda \deg(g)/\deg(\phi).$$ By Proposition \ref{propHensel} and Lemma \ref{lem:compareValuations}, $v_0^h(G_i-G_i^{(n)})>\sigma$, where $G_i$ is given by Proposition \ref{prop:factorization}. Since $\sigma > \delta_{\max}(g)\ge  \delta_0(G_i^{(n)})$, we deduce by induction from Corollary \ref{cor:precision} that the algorithm will recursively detect and compute all irreducible factors of $g$ with the suitable precision. The last statement is obvious.
\end{proof}

\e 
For $v$ of arbitrary rank, the algorithm will return a correct answer as soon as all involved $\gamma$'s satisfy 
$$
\delta_{\max}(g)\le \sup(m\gamma, m\in \N)
$$
since we can then compute a suitable integer $n\in \N$ at each call. This might be a weaker condition than in Corollary \ref{cor:converge}. However,  it's not clear that the approximants $G^{(n)}$ of $G\in\fff$ converge to $G$.

\e
The following corollary is immediate:

\begin{corollary}\label{cor:OkutsuEquivalent}
Suppose that $v$ has rank one. Let $g\in \oo[x]$ monic square-free with $\chr(k)=0$ or $\chr(k)>\deg(g)$.
If $f\in \oo[x]$ is monic and satisfies $v_0(g-f)>\delta_{\max}(g)$, then $g$ and $f$ have same OM-factorization.
\end{corollary}

\nn{\bf Remark. }Since we do not know a priori the bounds $\delta_0(g)$ and $\delta_{\max}(g)$, we start in practice with a small precision and check if it is sufficient to detect the right end-slope of the current Newton polygon (see e.g. \cite{PW1}). If not, we double the precision and restart all computations.


\subsubsection{Complexity issues}\label{ssec:complexity} In the discrete rank-one case, Poteaux and Weimann carried out an accurate analysis of the complexity of the Irreducibility and Factorization algorithms \cite{PW2}.

The extension of this analysis to the more general algorithms discussed so far is a delicate task, which goes beyond the scope of this paper.


\begin{thebibliography}{}


\bibitem{Ab} S.S. Abhyankar, \textit{Irreducibility criterion for germs of analytic
functions of two complex variables}, Adv. Math. \textbf{35} (1989), 190--257.

\bibitem{AM} S.S. Abhyankar, T. Moh, \emph{Newton-Puiseux Expansion and Generalized Tschirnhausen Transformation}, J. Reine Angew. Math. \textbf{260} (1973), 47--83.

\bibitem{VT}M. Alberich-Carrami$\tilde{\mbox{n}}$ana,  J. Gu\`ardia, E. Nart, J. Ro\'e, \emph{Valuative trees of valued fields}, preprint arXiv:2107.09813v3 [math.AG].

\bibitem{BN}M. dos Santos Barnab\'e, J. Novacoski, \emph{Valuations on $\kx$ approaching a fixed irreducible polynomial}, J. Algebra {\bf 592} (2022), 100--117.


\bibitem{Bauch}J.-D. Bauch, \emph{Computation of integral bases}, J. Number Th. {\bf 165} (2016), 382–-407.



\bibitem{endler}O. Endler, \emph{Valuation Theory}, Universitex, Springer-Verlag Berlin Heidelberg, 1972.







\bibitem{Gathen} J.v.z. Gathen, G. J\"urgen, \emph{Modern Computer Algebra}, Cambridge University Press, 2013.



\bibitem{okutsu}
J. Gu\`{a}rdia, J.  Montes, E.  Nart, \emph{Okutsu invariants and Newton polygons}, Acta Arith. {\bf 145} (2010), 83--108.

\bibitem{bordeaux}J. Gu\`{a}rdia, J.  Montes, E.  Nart, \emph{Higher  Newton polygons in the computation of discriminants and prime ideal decomposition in number fields}, J. Th\'eor. Nombres Bordeaux {\bf 23} (2011), no. 3, 667--696.

\bibitem{HN}J. Gu\`{a}rdia, J. Montes, E. Nart, \emph{Newton polygons of higher order in algebraic number theory}, Trans. Amer. Math. Soc.  {\bf 364} (2012), no. 1, 361--416.


\bibitem{newapp}J. Gu\`{a}rdia, J. Montes, E.  Nart, \emph{A new computational approach to ideal theory in number fields},   Found. Comput. Math. {\bf 13} (2013), 729--762.


\bibitem{bases}J. Gu\`{a}rdia, J. Montes, E.  Nart, \emph{Higher Newton polygons and integral bases}, J. Number Theory {\bf 147} (2015), 549-–589.


\bibitem{gen} J. Gu\`{a}rdia, E.  Nart, \emph{Genetics of polynomials over local fields}, in \emph{Arithmetic, geometry, and coding theory}, Contemp. Math. vol. 637 (2015), 207-241.

\bibitem{hos}F.J. Herrera Govantes, M.A. Olalla Acosta, M. Spivakovsky, \emph{Valuations in algebraic field extensions}, J. Algebra {\bf 312} (2007), no. 2, 1033--1074.

\bibitem{hmos}F.J. Herrera Govantes, W. Mahboub, M.A. Olalla Acosta, M. Spivakovsky, \emph{Key polynomials for simple extensions of valued fields}, preprint, arXiv:1406.0657v4 [math.AG], 2018.




\bibitem{Kuhl}F.-V. Kuhlmann, \emph{Value groups, residue fields, and bad places of rational function fields}, Trans. Amer. Math. Soc. {\bf 356} (2004), no. 11, 4559--4660.




\bibitem{mcla} S. Maclane, \emph{A construction for absolute values in polynomial rings}, Trans. Amer. Math. Soc. {\bf40} (1936), pp. 363--395.

\bibitem{mclb} S. Maclane, \emph{A construction for prime ideals as absolute values of an algebraic field}, Duke Mathematical Journal {\bf2} (1936), pp. 492--510.

\bibitem{montes}J. Montes, \emph{Pol\'{\i}gonos de Newton de orden superior y aplicaciones aritm\'eticas}, PhD Thesis, Universitat de Barcelona, 1999.



\bibitem{NMO} N. Moraes de Oliveira,
\emph{Inductive valuations and defectless polynomials over henselian fields}, PhD Thesis, Universitat Aut\`onoma de Barcelona, 2019.

\bibitem{defless} N. Moraes de Oliveira, E. Nart, \emph{Defectless polynomials over henselian fields and inductive valuations}, J. Algebra, {\bf 541} (2020), 270--307.

\bibitem{RPO} N. Moraes de Oliveira, E. Nart, \emph{Computation of residual polynomial operators of inductive valuations}, JPAA {\bf 225-9} (2021), 106668.



\bibitem{KP} E. Nart, \emph{Key polynomials over valued fields}, Publ. Mat. {\bf 64} (2020), 195--232.

\bibitem{MLV} E. Nart, \emph{Maclane-Vaqui\'e chains of valuations on a polynomial ring}, Pacific J. Math. {\bf 311-1} (2021), 165--195.

\bibitem{Rig} E. Nart, \emph{Rigidity of valuative trees under henselization}, arXiv:2202.0204v1 [math.AG], to appear in Pacific J. Math.

\bibitem{NN} E. Nart, J. Novacoski, \emph{The defect formula}, in preparation.




\bibitem{N2021} J. Novacoski, \emph{On Maclane--Vaqui\'e key polynomials}, J. Pure Appl. Algebra {\bf 225} (2021), 106644.

\bibitem{Ok}
K. Okutsu, \emph{Construction of integral basis I, II}, Proc. Japan Acad. Ser. A {\bf 58} (1982), 47--49, 87--89.

\bibitem{ore1} \O{}. Ore, \emph{Zur Theorie der algebraischen K\"orper}, Acta Math. {\bf44} (1923), pp. 219--314.

\bibitem{ore2} \O{}. Ore, \emph{Newtonsche Polygone in der Theorie der algebraischen K\"orper}, Math. Ann. {\bf 99} (1928), 84--117.




\bibitem{Pop} P. Popescu-Pampu, \emph{Approximate roots}, Fields Inst. Comm. {\bf 33} (2002), 1--37.

\bibitem{PW1} A. Poteaux, M. Weimann, \emph{A quasi-linear irreducibility test in $\mathbb{K}[[x]][y]$}, J. Comput. Comp. {\bf 31} (2022), no. 6, 1--52.

\bibitem{PW2} A. Poteaux, M. Weimann, \emph{Local polynomial factorisation: improving the Montes algorithm}, Proceedings of the 2022 ACM on International Symposium on
  Symbolic and Algebraic Computation ISSAC'22 (2022), 149--158.



\bibitem{hayden} H. D. Stainsby, \emph{Triangular bases of integral closures}, J. Symb. Comp. {\bf 87} (2018) 140--175.

\bibitem{Vaq0}M. Vaqui\'e, \emph{Famille admisse associ\'ee \`a une  valuation de $\kx$}, Singularit\'es Franco-Japonaises, S\'eminaires et Congr\'es 10, SMF, Paris (2005), Actes du colloque franco-japonais, juillet 2002, \'edit\'e par Jean-Paul Brasselet et Tatsuo Suwa, 391--428.

\bibitem{Vaq}
M. Vaqui\'e, \emph{Extension d'une valuation}, Trans. Amer. Math. Soc.  {\bf 359} (2007), no. 7, 3439--3481.

\bibitem{Vaq2}M. Vaqui\'e, \emph{Famille essential de valuations et d\'efaut d'une extension}, J. Algebra {\bf 311} (2007), no. 2, 859--876.

\bibitem{Vaq3}M. Vaqui\'e, \emph{Extensions de valuation et polygone de Newton}, Annales de l'Institute Fourier (Grenoble) {\bf 58} (2008), no. 7, 2503--2541.



\end{thebibliography}
\end{document}